\documentclass[12pt,letterpaper,twocolumn,english]{amsart}
\usepackage[margin=0.5in]{geometry}
\usepackage{amsmath}
\usepackage{amssymb}
\usepackage{amscd}

\usepackage[utf8]{inputenc}
\usepackage{graphicx}
\usepackage{amsfonts}
\usepackage{amstext}
\usepackage{mathrsfs}
\usepackage{bm}
\usepackage{tikz}
\usepackage{tikz-cd}
\usepackage[all]{xy}

\usepackage{amscd}
\usetikzlibrary{matrix,arrows,decorations.pathmorphing}
\tikzset{my loop/.style =  {to path={
  \pgfextra{}
  [looseness=12,min distance=10mm]
  \tikz@to@curve@path},font=\sffamily\small
  }}  
\usepackage{tikz-cd}

\usepackage{amscd}
\usetikzlibrary{matrix,arrows,decorations.pathmorphing}
\tikzset{my loop/.style =  {to path={
  \pgfextra{}
  [looseness=12,min distance=10mm]
  \tikz@to@curve@path},font=\sffamily\small
  }}  
\usepackage{tikz-cd}

\usepackage{enumitem}
\usetikzlibrary{arrows, decorations.pathmorphing}
\usepackage{hyperref}
\usepackage{comment}

\DeclareMathOperator{\SL}{SL}

\DeclareMathOperator{\Spec}{Spec}

\DeclareMathOperator{\tor}{tor}

\DeclareMathOperator{\GIT}{GIT}

\DeclareMathOperator{\ksba}{\mathrm{KSBA}}
\DeclareMathOperator{\ord}{\mathrm{ord}}
\DeclareMathOperator{\KE}{\mathrm{KE}}
\DeclareMathOperator{\Bl}{Bl}

\def\PP{\mathbb{P}}
\def\ZZ{\mathbb{Z}}
\def\RR{\mathbb{R}}
\def\CC{\mathbb{C}}
\def\QQ{\mathbb{Q}}
\def\BB{\mathbb{B}}
\def\LL{\mathbb{L}}
\def\hp{\mathfrak{H}}
\def\G{\Gamma}
\def\M{\mathbf{M}}
\def\co{\mathcal{O}}
\def\utau{\underline{\tau}}
\def\uo{\underline{0}}
\def\cx{\mathcal{X}}
\def\hm{\mathsf{h}}
\def\ua{\underline{\alpha}}
\def\btau{\bm{\tau}}
\def\ula{\underline{\lambda}}
\def\A{\mathcal{A}}
\def\vf{\varphi}
\def\uh{\underline{h}}
\def\uz{\underline{z}}
\def\D{\mathcal{D}}
\def\sm{\mathsf{m}}
\def\W{\mathcal{W}}
\def\fg{\mathfrak{g}}

\title{
Algebraic and analytic compactifications of moduli spaces
}
\author{Patricio Gallardo}
\address{
{\small Department of Mathematics,
University of California, Riverside,
900 University Ave.
Riverside, CA 92521
Skye Hall}
}
\email{pgallard@ucr.edu}

\author{Matt Kerr}
\address{{\small Department of Mathematics and Statistics
Washington University in St. Louis
Campus Box 1146
One Brookings Drive}}
\email{matkerr@wustl.edu}

\begin{document}
\sloppy 
\maketitle
\begin{abstract}
In this expository note, we offer an overview of the relationship between Hodge-theoretic and geometric compactifications of moduli spaces of algebraic varieties.
\end{abstract}

\section*{The tour ahead}

When studying algebraic objects, like subvarieties of a projective space, one notices that they are defined by polynomials, whose coefficients can vary.   This observation yields questions such as: 
\begin{itemize}[leftmargin=0.5cm]
\item When are objects with distinct coefficients equivalent?
\item What types of geometric objects appear if those coefficients move ``towards infinity''?  
\item Can we  make sense of ``equivalence classes at $\infty$''?
\end{itemize}
Searching for answers leads to the discovery of the so-called  \emph{moduli spaces} which parametrize equivalence classes of objects as diverse as varieties, vector bundles, representations of algebras, and so forth.

The construction of compact moduli spaces and the study of their geometry amount to a busy and central neighborhood of algebraic geometry nowadays.  Any vibrant district in an old city, of course, has too many landmarks to visit, and the first job of a tour guide is to curate a selection of sites and routes -- including multiple routes to the same site for the different perspectives they afford.

So we begin our tour along the unswerving path that compactifies algebro-geometric moduli spaces with ``limiting'' algebro-geometric invariants; the way is straight, but involves scaling a brick wall to discover what is meant by ``limits''.  Our subsequent turn down 19th-century vennels will unveil a connection as old as algebraic geometry itself:  associated to complex algebraic varieties are certain analytic invariants, called \emph{Hodge structures}, which yield a \emph{period map} from (the complex analytization of) our moduli space to a classifying space (for Hodge structures).  In the nicest cases, these latter spaces are locally symmetric varieties, the period map is either an isomorphism or an open embedding, and automorphic forms both invert the map and compactify the moduli space.  While lacking the ideological consistency of the former route, the ``limits'' of this one are more conceptually straightforward --- with Calculus providing a door in the wall.

For each given moduli problem there is an overwhelming wealth of routes to a compactification.  On the tour ahead, the reader will encounter various geometric and Hodge-theoretic methods that produce GIT, KSBA, $K$-stability, Baily-Borel, and toroidal compactifications of the same moduli space.  To travel towards understanding their differences --- and especially their spectacular coincidences --- is the aim of our brief journey.

\section{
Selected examples
}
\label{sec:ExampleModuliTheory}

Before discussing generalities any further, we turn to a series of examples (\S\S\ref{sec:ellipticCurves}-\ref{sec:Cubic4folds}) in which the period map is just given by a vector of complex (multivalued) functions.  
We will use each example to illustrate techniques and methods in moduli theory which will be discussed in the subsequent sections.  The unavoidable omission of many important cases and methods is a function of curating the tour, and we apologize to the experts in advance.

\subsection{ Elliptic curves} 
\label{sec:ellipticCurves}
With a rich history going back to Abel, Jacobi and Weierstra{\ss} in their guise as complex 1-tori, these are central objects in many areas of mathematics, from cryptography to complex analysis.  We'll use their moduli space (as cubic curves) to illustrate constructions such as geometric quotients, modular forms, and the period map.

We begin from an algebro-geometric perspective, with plane cubics --- that is, smooth algebraic curves defined as the complex solutions of homogeneous polynomials of degree three in 3 variables:
$$
F(x_0,x_1,x_2) = 
a_{300}x_0^3  + \ldots + 
a_{012}x_1x_2^2 + a_{003}x_2^3.
$$
Without the smoothness requirement, the coefficients of such equations comprise all ordered 10-tuples of complex numbers $[a_{300}: a_{210}: \cdots : a_{012}: a_{003}]$, not all zero, and defined up to scaling (by $\CC^{\times}$).  These are the complex points of the projective space $\PP^{9}$. 

If we are only interested in keeping track of the 
complex solutions of polynomial equations within projective space, then there is a solution to our classification problem: \emph{the Hilbert scheme}.  In this case, we need to fix an invariant known as the Hilbert polynomial which records geometric information such as their dimension and degree. 
It was shown in 1961 by Grothendieck that there exists a projective scheme $Hilb^{p(m)}_r$ which parametrizes all the closed 
complex solutions of polynomial equations in $\PP^r$ with Hilbert polynomial $p(m)$.
In our particular example, the Hilbert polynomial of plane cubics in $\PP^2$ is equal to $p(m)=3m$, and the associated Hilbert scheme is $\PP^9$.

However, the Hilbert scheme is not the solution that we are looking for. The reason is that given an elliptic curve defined by the equation $\{ F(x_0,x_1,x_2) = 0 \}$, we can use a linear change of coordinates
$x_i \mapsto a_{i0}x_0 + a_{i1}x_1 + a_{i2}x_2$ with      
$a_{ij} \in \mathbb{C}$, 
to obtain another equation $\{G(x_0, x_1, x_2) = 0 \}$. The elliptic curve defined by this second equation is isomorphic to the first one, and yet the Hilbert scheme tells us that they are different objects.  

The critical point here is that if we are parametrizing varieties $X \subset \PP^n$ with a fixed Hilbert polynomial, then we want to  account for the automorphisms of the ambient projective space. 
In our example, the ambient space of elliptic curves
$\{F(x_0, x_1, x_2) =0 \}$ is $\PP^2$, and  
the group associated to linear transformations 
$x_i \mapsto a_{i0}x_0 + a_{i1}x_1 + a_{i2}x_2$
has a dimension equal to 8. Therefore, of the 9 degrees of freedom associated with the coefficients of the equation $F(x_0,x_1,x_2)$, only  \emph{one} is intrinsic to the geometry of the elliptic curve, while the other eight are related to the linear change of coordinates in the ambient space.    

We arrive at our first (and almost correct!) definition of a moduli space:
let $\mathcal{U} \subset \mathbb{P}^{9}$ be the open locus parametrizing smooth elliptic curves, and $\SL(3, \CC)$ the group associated to linear change of coordinates among the variables
$x_0$, $x_1$ and $x_2$.  The quotient 
$$
\M_{1,1} := \mathcal{U}/\SL(3, \CC) \cong \mathbb{C}
$$
is the ``moduli space" of elliptic curves up to isomorphism.  
It will be tempting to consider a naive quotient $\PP^9/\SL_3(\mathbb{C})$ for constructing a compactification of $\M_{1,1}$. However, these naive quotients are usually either  of the wrong dimension or yield a non-Hausdorff topological space. 

The correct framework for constructing quotients within algebraic geometry is given by \emph{Geometric Invariant Theory (GIT)}, initiated by Mumford in 1969 \cite{mumford1994geometric}. A key result from GIT is the existence of a larger open locus   $\mathcal{U} \subset \mathcal{U}^{ss} \subset \PP^9$  and a well-defined quotient\footnote{ $U/G$ denotes a geometric quotient while $U/ \! \!/G$ denotes a categorical one. }
such that
$$ 
\overline{\M}_{1,1}^{\tiny\text{GIT}}
:=
\mathcal{U}^{ss}/ \! \!/ \SL_3(\mathbb{C}) 
\cong \PP^1
$$
is a complex projective variety compactifying $\M_{1,1}$. 
Moreover,  there is a unique ``minimal" closed $\SL_3(\mathbb{C})$-orbit associated to the point  $\overline{\M}_{1,1}^{\tiny\text{GIT}} \setminus \M_{1,1}$: namely, the orbit of the plane cubic
$
\{ x_0x_1x_2=0 \}.
$

Turning to a complex-analytic perspective, we can view a nonsingular cubic $C:=\{F(X_0,X_1,X_2)=0\}\subset \PP^2$ as a compact Riemann surface, with homology basis $\alpha,\beta\in H_1(C,\ZZ)$ (oriented so that $\alpha\cdot \beta=1$). Up to scale, there is a unique holomorphic form $\omega\in \Omega^1(C)$, with \emph{period ratio} $\tau:={\int_{\beta}\omega}/{\int_{\alpha}\omega}$ in the upper-half plane $\hp$.  This $\tau$ is well-defined modulo the action of $\gamma=\left(\begin{smallmatrix}a&b\\c&d\end{smallmatrix}\right)\in \G=\SL_2(\ZZ)$ through fractional-linear transformations $\gamma(\tau)=\tfrac{a\tau+b}{c\tau+d}$ induced by changing the homology basis.  We claim that the resulting (analytic) invariant $[\tau]\in \hp/\G$ captures the (algebraic) isomorphism class of $C$.

This is closely related to the classical theory of modular forms.  The (biperiodic) Weierstra{\ss} $\wp$-function associated to the lattice $\Lambda=\ZZ\langle 1,\tau\rangle$, 
$$\wp(u):=u^{-2}+\textstyle{\sum_{\lambda\in\Lambda\setminus 0}}[(u+\lambda)^{-2}-\lambda^{-2}],$$ 
satisfies $(\wp')^2=4\wp^3-g_2(\tau)\wp-g_3(\tau)$, where $g_2(\tau):=60\sum_{\lambda\in \Lambda\setminus 0}\lambda^{-4}$ and $g_3(\tau):=140\sum_{\lambda\in\Lambda\setminus 0}\lambda^{-6}$ are modular forms $M_k(\Gamma)$ of weights $k=4$ resp. $6$.  That is, they transform by the automorphy factor $(c\tau+d)^k$ under pullback by $\gamma\in \Gamma$, which makes $\bm{\jmath}:= \tfrac{g_2^3}{g_2^3-27g_3^2}\colon \hp/\Gamma\to \CC$ into a well-defined function.  Evidently, the image $E_{\tau}$ of the map $W\colon \CC/\Lambda \hookrightarrow \PP^2$ sending $u\mapsto [1:\wp(u):\wp'(u)]$ is a Weierstra{\ss} cubic 
\begin{equation}\label{I1}\small
X_2^2X_0=4X_1^3-g_2 X_1 X_0^2-g_3 X_0^3=4\textstyle{\prod_{i=1}^3}(X_1-e_i X_0),
\end{equation}
with period ratio $\tau$.

The key point is that any $C$ can be brought into this form (without changing $[\tau]$) through the action of $\SL_3(\CC)$ on coordinates. Fix a flex point $o\in C$ (i.e. $(C\cdot T_o C)=3$); since the dual curve $\check{C}$ has degree $6$, there are 3 more tangent lines $\{T_{p_i}C\}_{i=1}^3$ passing through $o$.  The $\{p_i\}$ are collinear, since otherwise one could construct a degree-1 map $C\to \PP^1$.  So we may choose  coordinates to have $o=[0:0:1]$, $T_o C=\{X_0=0\}$, $\{X_2=0\}\cdot C=\sum_{i=1}^3 p_i$, and $\sum_{i=1}^3\tfrac{X_1}{X_0}(p_i)=0$, which puts us in the above form \eqref{I1}.  In fact, if $\jmath:=\bm{\jmath}(\tau)\notin \{0,1,\infty\}=:\Sigma$, then rescaling yields a member of the family
\begin{equation}\label{I2}
y^2=4x^3-\tfrac{27\jmath}{\jmath-1}x-\tfrac{27\jmath}{\jmath-1}
\end{equation}
over $\PP^1\setminus \Sigma$, whose period map $[\tau]\colon 
\PP^1\setminus \Sigma \to \hp/\Gamma$ composed with 
$\bm{\jmath}$ extends to the identity $\CC\to\CC$.  Hence 
$C\underset{\text{\tiny SL}_3}{\cong}C'$ $\implies$ 
$[\tau]=[\tau']$ $\implies$ $\jmath=\jmath'$ $\implies$ 
$C\underset{\text{\tiny SL}_3}{\cong}C'$ yields the claimed 
equivalence of analytic and algebraic ``moduli'', and 
$\bm{\jmath}$ is an isomorphism.

While the choice of $o$ does not refine the moduli problem, keeping track of the ordered 2-torsion subgroup $\{o,p_1,p_2,p_3\}$ does.  In \eqref{I1}, this preserves the ordering of the $\{e_i\}$, which are parametrized by the weight-2 modular forms $\wp(\tfrac{1}{2}),\wp(\tfrac{\tau}{2}),\wp(\tfrac{\tau+1}{2})$ with respect to $\Gamma(2):=\ker\{\Gamma\to\mathrm{SL}_2(\ZZ/2\ZZ)\}$.  The roles of $\bm{\jmath}$ and \eqref{I2} are played by $\bm{\ell}:=\tfrac{e_3-e_2}{e_1-e_2}\colon \hp/\Gamma(2)\overset{\cong}{\to}\PP^1\setminus \Sigma$ and the Legendre family $y^2=x(x-1)(x-\ell)$, with $\bm{\jmath}=\tfrac{4}{27}\tfrac{(1-\bm{\ell}-\bm{\ell}^2)^3}{\bm{\ell}^2 (1-\bm{\ell})^2}$ describing the 6:1 covering $\hp/\Gamma(2)\twoheadrightarrow \hp/\Gamma$.  Notice that $\ell$ parametrizes the cross-ratio of 4 ordered points on $\PP^1$.

For any $N\geq 3$ we can let $\Gamma(N)\ltimes\ZZ^2$ act on $\hp\times \CC$ by $(\gamma,\lambda).(\tau,z):=(\gamma(\tau),\tfrac{z+\lambda}{c\tau+d})$ and take the quotient to produce the universal elliptic curve $\mathcal{E}(N)$ with level-$N$ structure (marked $N$-torsion) over the modular curve $Y(N):=\hp/\Gamma(N)$.  To produce an algebraic realization, we can use Jacobi resp. modular forms $M_k(\Gamma(N))$ to embed then in a suitable projective space.  (Indeed, $[g_2{:}g_3]\in \mathbb{WP}[4{:}6]$ and $[e_1{:}e_2{:}e_3]\in \PP^2$ already did this for $N=1$ and $2$.)  The compactification $\overline{Y}(N)$ of $Y(N)$ so obtained adds $\tfrac{N^2}{2}\prod^{\text{$p$ prime}}_{p\mid N}(1-\tfrac{1}{p^2})$ points called \emph{cusps}, over which the elliptic fiber degenerates to a cycle of $N$ $\PP^1$'s; in fact, we have $
\overline{Y}(N) \setminus Y(N)= \PP^1(\QQ)/\Gamma(N)$.  Going around a cusp subjects a basis of integral homology to a transformation conjugate to $\left(\begin{smallmatrix}1&N\\0&1\end{smallmatrix}\right)$.

Now consider an algebraic realization $\mathcal{C}\to \M$ of $\mathcal{E}(N)\to Y(N)$; e.g. for $N=3$, the Hesse pencil $tX_0X_1X_2=X_0^3+X_1^3+X_2^3$ over $(t\in )\CC\setminus\{1,\zeta_3,\bar{\zeta}_3\}$ has a marked $(\ZZ/3\ZZ)^2$ subgroup as base-locus (where the curves meet the coordinate axes).  The \emph{monodromy group} generated by \emph{all} loops in $\M$ (acting on $H_1$ of some fiber) is tautologically $\Gamma(N)$.  So the period ratio $\tau$ yields a well-defined \emph{period map} $\M \to Y(N)$ inverted by modular forms (as for $N=1$ and $2$), exchanging algebraic and analytic moduli.  The other moral here is that refining the moduli problem (e.g. level structure) produces smaller monodromy group, hence more boundary components (in this case, cusps) in the compactification.

\subsection{Picard curves and points in a line}
\label{sec:pointsP1}
Next, we'll visit ideas that are central for constructing compact moduli spaces. They include the use of finite covers to associate periods to varieties without periods and  the first example of ``stable pairs." 

We begin this time with the complex-analytic point of view.  Though ordered collections of $n$ points in $\PP^1$ do not themselves have periods, we can consider covers $C\twoheadrightarrow \PP^1$ branched over such collections, generalizing the $n=4$ case of Legendre elliptic curves.  For $n=5$, compactifying
\begin{equation}\label{I3}
y^3=x^4+G_2 x^2 +G_3 x+G_4=\textstyle{\prod_{i=1}^4} (x-t_i)
\end{equation}
to $C\subset \PP^2$ yields a genus-3 curve with cubic automorphism 
$\mu\colon y\mapsto \zeta_3 y$ and (up to scale) unique holomorphic 
differential $\omega=\tfrac{dx}{y}$ with $\mu^*\omega=\bar{\zeta}_3 \omega$. 
The moduli spaces 
$\M_{\tiny\text{ord}}:=\PP(\sum t_i=0)\setminus 
\cup_{i<j}\{t_i=t_j\}\cong \PP^2\setminus \{\text{6 lines}\}$ and 
$\M=\M_{\tiny\text{ord}}/\mathfrak{S}_4$ parametrize ordered 5-tuples in $\PP^1$ (via $\{t_i\}$) resp. unordered 4-tuples in $\CC$ (via $\{G_j\}$).

To define period maps, we first describe the monodromy group $\Gamma$ through which $\pi_1(\M,\mathsf{m})$ acts on $H_1(C_{\mathsf{m}},\ZZ)$. As it sends symplectic bases to symplectic bases and must be compatible with $\mu$, it should be plausible that this takes the form
$$\Gamma=\mathrm{Sp}_6(\ZZ)^{[\mu_*]}\cong U((2,1);\co),\;\;\;\;\;\co=\ZZ\langle 1,\zeta_3\rangle;$$
and that for $\M_{\tiny\text{ord}}$ this is replaced by the subgroup $\Gamma(\sqrt{-3})=U((2,1);\sqrt{-3}\co)$.  Now given a basis $\alpha,\beta_1,\beta_2\in H_1(C,\ZZ)_{\mu_*=\bar{\zeta}_3}$, the period vector $\utau:=(\tau_1,\tau_2)$ ($\tau_i=\int_{\beta_i}\omega/\int_{\alpha}\omega$) lies in a 2-ball $\BB_2$ (by Riemann's $2^{\text{nd}}$ bilinear relation), on which $\gamma=(\gamma_{ij})_{i,j=0}^2$ acts via $\gamma(\utau):=$
$$\j_{\gamma}(\utau)^{-1}.\left(\gamma_{10}+\gamma_{11}\tau_1+\gamma_{12}\tau_2,\,\gamma_{20}+\gamma_{21}\tau_1+\gamma_{22}\tau_2\right)$$
with $\j_{\gamma}(\utau):=\gamma_{00}+\gamma_{01}\tau_1+\gamma_{02}\tau_2$.  So we get period maps $\phi\colon\M\to\BB_2/\Gamma$ and $\tilde{\phi}\colon \M_{\tiny\text{ord}} \to\BB_2/\Gamma(\sqrt{-3})$, whose images omit 1 resp. 6 disk-quotients.  Writing $\M_{\tiny\text{ord}}'$ resp. $\M'$ for $\PP(\sum t_i=0)\setminus\cup_{i<j<k}\{t_i=t_j=t_k\}\cong \PP^2\setminus\{\text{4 pts.}\}$ and its $\mathfrak{S}_4$-quotient, these maps extend to isomorphisms $\M_{\tiny\text{ord}}'\overset{\cong}{\to}\BB_2/\Gamma(\sqrt{-3})$ and $\M'\overset{\cong}{\to}\BB_2/\Gamma$.

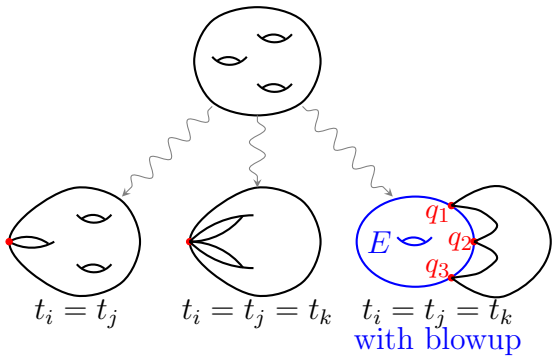
\begin{figure}[h!]
\centering
\begin{tikzpicture}[scale=0.6]
\draw [thick] plot [tension=0.6,smooth cycle] coordinates {(1.4,4) (1,4.9) (0,5.2) (-1,4.9) (-1.4,4) (-1,3.1) (0,2.8) (1,3.1)};
\draw [thick] (-0.9,4) arc (130:50:0.4);
\draw [thick] (-1,4.1) arc (-140:-40:0.5);
\draw [thick] (0.1,3.4) arc (130:50:0.4);
\draw [thick] (0,3.5) arc (-140:-40:0.5);
\draw [thick] (0.1,4.5) arc (130:50:0.4);
\draw [thick] (0,4.6) arc (-140:-40:0.5);
\draw [thick] plot [tension=0.6,smooth cycle] coordinates {(-2.6,0) (-3,0.9) (-4,1.2) (-5,0.7) (-5.5,0) (-5,-0.7) (-4,-1.2) (-3,-0.9)};
\draw [thick] (-5.5,0) arc (120:54:0.8);
\draw [thick] (-5.5,0) arc (-120:-60:1);
\draw [thick] (-3.9,-0.6) arc (130:50:0.4);
\draw [thick] (-4,-0.5) arc (-140:-40:0.5);
\draw [thick] (-3.9,0.5) arc (130:50:0.4);
\draw [thick] (-4,0.6) arc (-140:-40:0.5);
\filldraw [red] (-5.5,0) circle (2pt);
\node at (-4,-1.6) {$t_i=t_j$};
\draw [thick] plot [tension=0.6,smooth cycle] coordinates {(1.4,0) (1,0.9) (0,1.2) (-1,0.7) (-1.5,0) (-1,-0.7) (0,-1.2) (1,-0.9)};
\filldraw [red] (-1.5,0) circle (2pt);
\draw [thick] (-1.5,0) arc (135:90:2);
\draw [thick] (-1.5,0) arc (-135:-90:2);
\draw [thick] (-1.5,0) arc (-90:-38:1.5);
\draw [thick] (-1.5,0) arc (90:38:1.5);
\node at (0,-1.6) {$t_i=t_j=t_k$};
\draw [blue,thick] (3.5,0) ellipse (1.3cm and 1cm);
\draw [blue,thick] (3.2,0) arc (130:50:0.4);
\draw [blue,thick] (3.1,0.1) arc (-140:-40:0.5);
\filldraw [red] (4.3,0.8) circle (2pt);
\filldraw [red] (4.3,-0.8) circle (2pt);
\filldraw [red] (4.8,0) circle (2pt);
\draw [thick] plot [tension=0.8,smooth] coordinates {(4.3,0.8) (5.5,1.2) (6.5,0) (5.5,-1.2) (4.3,-0.8)};
\draw [thick] plot [smooth] coordinates {(4.3,0.8) (4.9,0.7) (5.3,0.4) (5.2,0.2) (4.8,0)};
\draw [thick] plot [smooth] coordinates {(4.3,-0.8) (4.9,-0.7) (5.3,-0.4) (5.2,-0.2) (4.8,0)};
\node at (4,-1.6) {$t_i=t_j=t_k$};
\node [blue] at (4,-2.2) {with blowup};
\node [red] at (4,0.6) {\small $q_1$};
\node [red] at (4.5,0) {\small $q_2$};
\node [red] at (4,-0.6) {\small $q_3$};
\draw[gray,-stealth,decorate,decoration={snake,amplitude=3pt,pre length=2pt,post length=3pt}] (-1,3) -- (-3,1);
\draw[gray,-stealth,decorate,decoration={snake,amplitude=3pt,pre length=2pt,post length=3pt}] (0,2.8) -- (0,1.2);
\draw[gray,-stealth,decorate,decoration={snake,amplitude=3pt,pre length=2pt,post length=3pt}] (1,3) -- (3,1);
\node [blue] at (2.7,0) {$E$};
\end{tikzpicture}
\caption{Degenerating a Picard curve}
\label{fig:Pic}
\end{figure}

The Picard modular forms\footnote{Here $M_k(\BB_2,\Gamma_0)$ comprises holomorphic functions on $\BB_2$ satisfying $f(\gamma(\utau))=\j_{\gamma}(\utau)^k f(\utau)$ for all $\gamma\in \Gamma_0$; and $S\Gamma_0:=\ker(\det)\cap \Gamma_0$.} describing their inverses are none other than the $t_i\in M_3(\BB_2,\mathrm{S}\Gamma(\sqrt{-3}))$ and $G_j\in M_{3j}(\BB_2,\mathrm{S}\Gamma)$ in \eqref{I3}.  The compactifications of the ball quotients so obtained add only (4 resp. 1) points, extending (say) $\tilde{\phi}$ to $\tilde{\phi}^*\colon \PP^2\overset{\cong}{\to}(\BB_2/\Gamma(\sqrt{-3}))^*$.  To understand the meaning of $\tilde{\phi}^*$, notice that colliding two $t_i$ in \eqref{I3} and normalizing yields a genus 2 curve with cubic automorphism, whose (single) period ratio is parametrized by one of the disk-quotients previously omitted.  When 3 $t_i$ collide in \eqref{I3}, the normalization has genus 0 and thus no moduli, which explains the 4 boundary points in $(\BB_2/\Gamma(\sqrt{-3}))^*$.  The unnormalized scenarios are depicted in the left and middle degenerations in Fig. \ref{fig:Pic}.

But this is not the only way to approach the collision of 3 $t_i$.  After a linear change in coordinates, the (2-parameter) degeneration takes the form
$$y^3=(x-s_1)(x-s_2)(x+s_1+s_2)(x-1)$$
in a neighborhood of $\underline{s}=\uo$. Restricting to $s_i=a_i t$ ($|t|<\epsilon$, $a_i\in \CC$ fixed) yields a 1-parameter family $\cx\to \Delta$ over a disk.  Blowing $\cx$ up at $(x,y,t)=\uo$ produces the exceptional divisor
$$E\colon\;\;Y^3=(X-a_1Z)(X-a_2Z)(X+(a_1+a_2)Z),$$
which is an elliptic curve with period ratio $\zeta_3$.  The singular fiber of the blowup is the union of $E$ with the normalization ($\cong \PP^1$) of $y^3=x^3(x-1)$, glued along $E\cap\{Z=0\}=\{q_1,q_2,q_3\}$ (see the rightmost degeneration in Fig. \ref{fig:Pic}), with $\mu$ acting on the lot (and cyclically permuting the $q_i$).  While the modulus of $E$ is just $[\zeta_3]$, the ratios $\eta_i$ of the semiperiods $\int_{q_i}^{q_{i+1}}\omega_E$ to a period of $E$ vary in $[a_1{:}a_2]\in \PP^1$, and are related by complex multiplication by $\zeta_3$ (i.e. $\mu|_E$).  In fact, as $t\to 0$ it turns out that (for some choice of $\{\alpha,\beta_1,\beta_2\}$) $\tau_1\sim \tfrac{\log(t)}{2\pi\sqrt{-1}}$ blows up, while $\tau_2$ limits to (say) $\eta_1(\underline{a})$, a limit which becomes well-defined in $E/\langle\mu\rangle \cong \PP^1$.

The upshot is that if we replace the 4 boundary points of $(\BB_2/\Gamma(\sqrt{-3}))^*$ by copies of $E/\langle \mu\rangle$, then these semiperiod ratios extend $\tilde{\phi}$ to an isomorphism from $\mathrm{Bl}_{\{\text{4 pts.}\}}(\PP^2)$ to the resulting $(\BB_2/\Gamma(\sqrt{-3}))^{**}$.  (We have to blow up at the 4 triple-intersection points to make $a_2/a_1$ well-defined.)  This is a first example of using limiting mixed Hodge structures (LMHS) to extend period maps to a \emph{toroidal} compactification ``$**$'' refining the Baily-Borel ``$*$'' compactification.  

This close relationship between moduli of $n=5$ points, Picard curves, and ball quotients has a beautiful extension to configurations of up to 12 points by Deligne and Mostow.  The geometric meaning of the corresponding ``*'' and ``**'' compactifications will be addressed in \S3.

Going back to the moduli of $n$ ordered points in $\PP^1$, and adopting a geometric viewpoint, we should phrase the problem in terms of objects up to an equivalence relation.  An ``object'' here is an $n$-pointed curve $(\PP^1,(p_1,\ldots,p_n))$, which is equivalent to $(\PP^1,(q_1,\ldots,q_n))$ if $g(p_i)=q_i$ ($1\leq i\leq n$) for some $g\in \mathrm{Aut}(\PP^1)$.  The resulting moduli space \begin{align*}\M_{0,n}&\phantom{:}=\textstyle{\left(\left( \PP^1 \right)^n \setminus \bigcup_{i < k} \Delta(ik)\right)}/  \SL_2(\mathbb{C})\\\Delta(ik) &:=\{ (x_1, \ldots, x_n) \in (\PP^1)^n\; | \; x_i{=} x_k  \;\; 1 \leq i {,}\; k \leq n\}\end{align*} is an $(n-3)$-dimensional quasi-projective variety.

It was discovered in the late 1960s by Grothendieck and later by Knudsen that we can define pairs of a more general sort, called
\emph{stable $n$-pointed curves of genus $0$}, whose moduli space $\overline{\M}_{0,n}$ is a smooth, compact algebraic variety.  Moreover, the complement $\overline{\M}_{0,n} \setminus \M_{0,n}$ is a normal crossing divisor with smooth irreducible components. 

This new type of ``stable pair" added at the boundary of the compactification is a connected, but possibly
reducible, complex curve $C$ together with $n$ smooth distinct labelled points $p_1, \ldots,  p_n$ in $C$, satisfying the following conditions:
\begin{itemize}[leftmargin=0.5cm]
\item The arithmetic genus of $C$ is equal to $0$.
\item $C$ has only ordinary double points and every irreducible component of $C$ is isomorphic to the projective line $\PP^1$.
\item On each component of $C$ there are at least three points which are either one of the marked points $p_i$ or a double point, i.e., the intersection of two components of $C$. 
\end{itemize}
Let's consider the case of $n=5$ closely.  The moduli space $\overline{\M}_{0,5}$ is two dimensional and isomorphic to the blow-up of 
$\PP^2$ at four points in general position.  The $10$ irreducible divisors $D_I$ in the boundary $\overline{\M}_{0,5} \setminus \M_{0,5}$ are labelled by subsets $I \subset \{1, 2,\ldots, 5\}$ with  $|I|=2$.
\begin{figure}[h!]
\centering
\begin{tikzpicture}[scale=0.6]
\def\ra{-3.5}
\def\rb{3.2}
\draw [thick] plot [tension=0.6,smooth cycle] coordinates {(1.4+\ra,4) (1+\ra ,4.9) (0+\ra,5.2) (-1+\ra,4.9) (-1.4+\ra,4) (-1+\ra,3.1) (0+\ra,2.8) (1+\ra,3.1)};
\node [red] at (0.5+\ra,4.8) {\small $1$};
\filldraw [red] (0.5+\ra,4.5) circle (2pt);
\node [red] at (-0.6+\ra,4.7) {\small $3$};
\filldraw [red] (-0.3+\ra,4.7) circle (2pt);
\node [red] at (0.7+\ra,3.8) {\small $2$};
\filldraw [red] (0.7+\ra,3.5) circle (2pt);
\node [red] at (-0.8 +\ra,3.9) {\small $4$};
\filldraw [red] (-0.5+\ra,3.9) circle (2pt);
\node [red] at (-0.2 +\ra,3.2) {\small $5$};
\filldraw [red] (-0.5+\ra,3.2) circle (2pt);

\draw [thick] plot [tension=0.6,smooth cycle] coordinates {(1.4+\rb,4) (1+\rb ,4.9) (0+\rb,5.2) (-1+\rb,4.9) (-1.4+\rb,4) (-1+\rb,3.1) (0+\rb,2.8) (1+\rb,3.1)};
\draw[gray,-stealth,decorate,decoration={snake,amplitude=3pt,pre length=2pt,post length=3pt}] (-1.8,4) -- (1.5,4);
\node [gray] at (0,3) {
};
\def\rc{6.0}
\draw [thick] plot [tension=0.6,smooth cycle] coordinates {(1.4+\rc,4) (1+\rc ,4.9) (0+\rc,5.2) (-1+\rc,4.9) (-1.4+\rc,4) (-1+\rc,3.1) (0+\rc,2.8) (1+\rc,3.1)};
\node [red] at (0.5+\rc,4.8) {\small $1$};
\filldraw [red] (0.5+\rc,4.5) circle (2pt);
\node [red] at (-0.3+\rb,4.7) {\small $3$};
\filldraw [red] (-0.0+\rb,4.7) circle (2pt);
\node [red] at (0.7+\rc,3.8) {\small $2$};
\filldraw [red] (0.7+\rc,3.5) circle (2pt);
\node [red] at (-0.8 +\rb,3.9) {\small $4$};
\filldraw [red] (-0.5+\rb,3.9) circle (2pt);
\node [red] at (-0.3 +\rb,3.2) {\small $5$};
\filldraw [red] (-0.0+\rb,3.2) circle (2pt);
\filldraw [blue] (1.4+\rb,4) circle (3pt);
\end{tikzpicture}
\caption{Generic limit  parametrized by $D_{12}$}
\label{fig:Mon}
\end{figure}
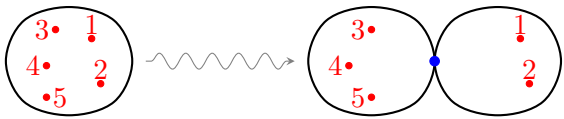
For example, the divisor $D_{12}$ generically parametrizes the union of two $\PP^1$s with the points distributed as in Fig. \ref{fig:Mon}. 

The reader may wonder if we can also compactify these $\M_{0,n}$ spaces by going down the same road as for elliptic curves. The answer is yes --- there are indeed GIT compactifications --- but with a new twist. 
We determined already that $\M_{0,n}$ is a quotient of an open locus within $( \PP^1 )^n$ by $\SL_{2}(\mathbb{C})$. Geometric Invariant Theory and subsequent developments imply that there are finitely many open loci $\mathcal{U}^{ss}_{\mathbf{w}}$, depending on a collection of rational numbers
$
\mathbf{w} = (w_1, \ldots, w_n)
$
with $0 < w_i \leq 1$ and $w_1 + \cdots + w_n=2$,
such that 
\begin{align*}
\textstyle{( \PP^1 )^n \setminus \bigcup_{i \leq k} \Delta(ik)
\subset 
\mathcal{U}^{ss}_{\mathbf{w}}
\subset 
( \mathbb{P}^1 )^n}
\end{align*}
and the quotient
\begin{align*}
\textstyle{( \PP^1 )^n / \! \!/_{\mathbf{w}} \SL_2 (\mathbb{C})
:=
\mathcal{U}^{ss}_{\mathbf{w}}
/ \! \!/ \SL_2\left( \mathbb{C} \right)}
\end{align*}
is a projective variety compactifying $\M_{0,n}$. For the case $n=5$ and depending on the choice of ``weights'' $\mathbf{w}$,  the quotients $( \PP^1 )^5 / \! \!/_{\mathbf{w}} \SL_2 (\mathbb{C})$ can be either $\PP^2$, $\PP^1 \times \PP^1$, or a blow-up of $\PP^2$ at $k$ points in general position with $1 \leq k \leq 4$.

\begin{figure}[h!]
\centering
\begin{tikzpicture}[scale=0.7]
\def\ra{-3.5}
\def\rb{3.2}
\draw [thick] plot [tension=0.6,smooth cycle] coordinates {(1.4+\ra,4) (1+\ra ,4.9) (0+\ra,5.2) (-1+\ra,4.9) (-1.4+\ra,4) (-1+\ra,3.1) (0+\ra,2.8) (1+\ra,3.1)};
\node [red] at (0.5+\ra,4.8) {\small $1$};
\filldraw [red] (0.5+\ra,4.5) circle (2pt);
\node [red] at (-0.6+\ra,4.7) {\small $3$};
\filldraw [red] (-0.3+\ra,4.7) circle (2pt);
\node [red] at (0.7+\ra,3.8) {\small $2$};
\filldraw [red] (0.7+\ra,3.5) circle (2pt);
\node [red] at (-0.8 +\ra,3.9) {\small $4$};
\filldraw [red] (-0.5+\ra,3.9) circle (2pt);
\node [red] at (-0.2 +\ra,3.2) {\small $5$};
\filldraw [red] (-0.5+\ra,3.2) circle (2pt);
\draw [thick] plot [tension=0.6,smooth cycle] coordinates {(1.4+\rb,4) (1+\rb ,4.9) (0+\rb,5.2) (-1+\rb,4.9) (-1.4+\rb,4) (-1+\rb,3.1) (0+\rb,2.8) (1+\rb,3.1)};
\draw[gray,-stealth,decorate,decoration={snake,amplitude=3pt,pre length=2pt,post length=3pt}] (-1.8,4) -- (1.5,4);
\node [gray] at (0,3) {
};
\def\rc{6.5}
\node [red] at (-0.3+\rb,4.7) {\small $3$};
\filldraw [red] (-0.0+\rb,4.7) circle (2pt);
\node [red] at (0.7+\rb,4.1) {\small $1=2$};
\filldraw [red] (0.7+\rb,3.8) circle (2pt);
\node [red] at (-0.8 +\rb,3.9) {\small $4$};
\filldraw [red] (-0.5+\rb,3.9) circle (2pt);
\node [red] at (-0.3 +\rb,3.2) {\small $5$};
\filldraw [red] (-0.0+\rb,3.2) circle (2pt);
\end{tikzpicture}
\caption{GIT limit  when $w_1 +w_2 < 1$ }
\label{fig:GITpoints}
\end{figure}
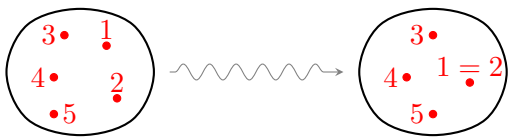

These GIT quotients are philosophically different from $\overline{\M}_{0,n}$, because they parametrize different types of geometric objects.
Remember that $\overline{\M}_{0,n}$ allows $\PP^1$ itself to degenerate, so as to keep the points distinct (as in the previous figure).  The GIT quotient, on the other hand, enables the points to collide amongst themselves in a controlled manner, and $\PP^1$ does not degenerate.  This scenario is depicted in  Fig. \ref{fig:GITpoints}.

And so we arrive at one of the main questions within moduli theory. \emph{How are} a priori \emph{different compactifications of a moduli space related to each other?} In our case, it is a theorem of Kapranov that for every $\mathbf{w}$ as above there is a morphism 
$\overline{\M}_{0,n} \to 
( \PP^1 )^n / \! \!/_{\mathbf{w}} \SL_2 (\mathbb{C})
$. 
A framework known as ``variation of GIT'' (VGIT), developed by Dolgachev, Hu, and Thaddeus, shows that changing the values of $\mathbf{w}$ induces birational transformations among the GIT quotients. For example, for each $n$ there are choices of $\mathbf{w}$ 
that yield $\PP^{n-3}$ and $(\PP^1)^{n-3}$ as quotients.

\subsection{Cubic surfaces}
\label{sec:cubicsurfaces}

In our first higher dimensional example, we will visit the idea of K-stability, which leads to a big part of current moduli theory.

Cubic surfaces $S_F \subset \PP^3$ are cut out by homogeneous polynomials $F(X_0,\ldots,X_3)$ of degree three.  The corresponding Hilbert scheme is equal to 
$\PP^{19}\cong\mathbb{P}H^0(\PP^3,\co(3))$ 
and the group $\SL_4( \mathbb{C})$ associated with the linear change of coordinates has dimension $15$.
Let $\mathbb{P}H^0(\PP^3,\co(3))^{\text{sm}}\subset \PP^{19}$
be the open locus parametrizing smooth cubic surfaces. 
The four-dimensional quotient 
\begin{align*}
\M:=
\mathbb{P}
H^0(\PP^3,\co(3))^{\text{sm}}
/
\mathrm{SL}_4(\CC)
\end{align*}
is our moduli space of smooth cubic surfaces. 

For compactifying $\M$, we observe that smooth cubic surfaces belong to a family known as \emph{Fano varieties}. These are complex projective varieties whose anticanonical bundle --- the determinant of the tangent bundle --- is ample.  A fundamental problem in complex geometry is to classify the Fano varieties that admit a K\"{a}hler-Einstein (KE) metric:  that is,  a K\"{a}hler metric $\omega$ with Ricci curvature
equal to the metric  ($\operatorname{Ric}(\omega) = \omega$).  
A celebrated result of Tian from 1990 implies that all \emph{smooth} cubic surfaces admit a KE metric.  So we arrive at a question which leads to much of the moduli theory for Fano varieties:  \emph{Can we construct a compactification 
$\M \subset \overline{\M}^{\KE}$ 
of $\M$ whose limiting objects admit a KE metric as well?}

The answer revolves around the relatively new concept of \emph{K-(semi)stability} first introduced in 1997 by Tian and then expanded in 2002 by Donaldson.  Leaving definitions aside for the moment, we only remark that $K$-(semi)stability induces substantial restrictions on a Fano variety's singularities. For example, if $X$ is a ``reasonable" Fano $K$-semistable variety, then $X$ should be normal.  

In our particular case, results of Odaka, Spotti and Sun from 2012 imply that the moduli of KE cubic surfaces is isomorphic to the GIT quotient:
 \begin{align*}
\overline{\M}^{\KE}
\cong
\overline{\M}^{\tiny\text{GIT}}
:=
\left( 
\mathbb{P}H^0(\PP^3,\co(3))
\right)^{\text{ss}} 
/ \! \! / \SL_4( \mathbb{C})
 \end{align*}
where  $\left( \mathbb{P}H^0(\PP^3,\co(3))\right)^{\text{ss}} \subset \PP H^0(\PP^3,\co(3))$
is the open subset parametrizing \emph{semistable} cubic surfaces: those which are either smooth, or admit $A_1$ or $A_2$ singularities (of the local form $x^2 +y^2 +z^2$ resp. $x^2+y^2+x^3$).
The compactification $\overline{\M}^{\tiny\text{GIT}}$ contains an open subset
$\M^{\text{s}}=\left( \mathbb{P}H^0(\PP^3,\co(3))\right)^{\text{s}} /\SL_4(\CC)$ called the \emph{stable locus}.
In our case, each point in the  stable locus corresponds to a unique isomorphism class of cubic surfaces with at worst $A_1$ singularities. 
Our GIT quotient compactifies the open stable locus by adding a single point $\omega$.  Although there is not a unique 
$\SL_4(\CC)$-orbit within $\PP^{19}$ associated to $\omega$, there is a unique closed one! In this case, the associated surface is the  cubic surface $\{ x_0x_1x_2 + x_3^3 =0 \}$ with three $A_2$ singularities. 

We should be proud of progress until now, and yet there is an aspect missing from the above story.  It is known since work of Cayley and Salmon circa 1850 that every smooth cubic surface has $27$ lines on it.  Can we produce a compactification that keeps tracks of such lines and their limits?  With that purpose, we choose in order $6$ disjoint lines out of the $27$.
A smooth cubic surface together with such a choice is called a marked cubic surface.
The resulting open moduli space $\M_{\tiny\text{ord}}$ has a finite map to $\M$ defined by forgetting the labeling of the lines.  The associated finite field extension
$\mathbb{C}\left( \M \right)    
\longrightarrow
\mathbb{C}\left( \M_{\tiny\text{ord}} \right)    
$
induces a unique compactificiation 
$\overline{\M}_{\tiny\text{ord}}^{\tiny\text{GIT}}$ such that
\begin{align*}
\xymatrix{ 
\M_{\tiny\text{ord}}
\ar[d]^{/W(E_6)} 
\ar@{^{(}->}[r]
& 
\overline{\M}_{\tiny\text{ord}}^{\tiny\text{GIT}}
\ar[d]^{/W(E_6)}
\\
\M  \ar@{^{(}->}[r]
& 
\overline{\M}^{\tiny\text{GIT}}
}    
&
\end{align*}
where we used the classical fact that the symmetric group of the 27 lines is the Weyl group of $E_6$.

Returning to a Hodge-theoretic perspective, we note that $S_F$ has no holomorphic 1- or 2-forms; so the situation is analogous to that of point-configurations on $\PP^1$, with no ``straightforward'' periods.  Recycling the cyclic-cover trick from \S\ref{sec:pointsP1} yields cubic threefolds
$$T_F:=\{X_4^3=F(X_0,\ldots,X_3)\}\subset \PP^4$$
with automorphism $\mu\colon X_4\mapsto \zeta_3 X_4$.  Up to scale, $T_F$ has a unique closed 3-form $\omega$ of type $(2,1)$ with $\mu^*\omega=\bar{\zeta}_3\omega$, whose periods over a basis of the $\co$-lattice $H_3(T_F,\ZZ)$ produce a point $[\utau]$ in a 4-ball $\BB_4$.

We now summarize some results of Allcock, Carlson and Toledo.  Note that the intersection form on $H_3(T_F,\ZZ)$ gives rise to a Hermitian form $\hm$ of signature $(4,1)$.  It turns out that the monodromy group (through which $\pi_1(\M)$ acts on $H_3$) is $\Gamma:=U(\hm,\co)$, the resulting period map $\phi\colon \M\to \BB_4/\Gamma$ is an embedding, and $\phi$ extends to an isomorphism $\M^{\text{s}}\cong \BB_4/\Gamma$ by allowing $A_1$ singularities in $S_F$.  Next, recalling that any smooth $S_F$ has 27 lines, we may refine our moduli (in analogy with the level-structures of \S\ref{sec:ellipticCurves}) by an ordered choice of 6 disjoint lines.  Writing $\M_{\tiny\text{ord}}$ ($\twoheadrightarrow \M$ by forgetting the lines) for the resulting moduli space, $\phi$ lifts to 
$\tilde{\phi}\colon \M_{\tiny\text{ord}}
\hookrightarrow 
\BB_4/\Gamma(\sqrt{-3})$, 
where $\Gamma(\sqrt{-3}):=U(\hm,\sqrt{-3}\co)$ and $\Gamma/\Gamma(\sqrt{-3})$ is again $W(E_6)$.  (Notice the analogy between the role of the marking and of the level structures in \S\ref{sec:ellipticCurves}.)   We postpone discussing the compactifications, but note that Allcock and Freitag have given an automorphic description of Cayley's cross-ratios (hence of $\tilde{\phi}^{-1}$).

\subsection{$K3$ surfaces}
\label{sec:K3surfaces}
Named in honor of Kummer, K\"ahler and Kodaira, these are simply connected compact complex 2-manifolds $X$ with $\Omega^2_X\cong \co_X$, hence (up to scale) a unique and nowhere-vanishing holomorphic 2-form $\omega\in \Omega^2(X)$.  Algebraic (quasi-polarized) $K3$s have countably many (19-dimensional) coarse moduli spaces $\M_g^{K3}$, one for every ``genus'' $g\geq 2$.   An open subset $\M_g^{K3,\circ}$ parametrizes surfaces $X\subset \PP^g$ of degree $2g-2$ for $g\geq 3$, and double covers $X\twoheadrightarrow \PP^2$ branched along a sextic curve for $g=2$.  Given a hyperplane section $h\subset X$ of one of these \emph{polarized} $K3$s, $V_{\ZZ}:=H_2(X,\ZZ)/[h]$ is a lattice of rank 21 (and signature $(2,19)$ under the intersection form), comprising classes of topological 2-cycles.  As $\int_h \omega=0$, integrating $\omega$ over (some basis of) these yields a well-defined period point $[\omega]\in \PP V_{\CC}^{\vee}\cong \PP^{20}$, which is further restricted by $\int_X\omega\wedge\omega=0$ (and $\int_X\omega\wedge\bar{\omega}>0$) to lie in a 19-dimensional submanifold $D$ called a type-IV (Hermitian symmetric) domain.  By a theorem of Piatetskii-Shapiro and Shafarevich, this period map $\phi$ embeds $\M_g^{K3,\circ }\hookrightarrow D/\Gamma_g$ (where $\Gamma_g=O(V_\ZZ)$) as the complement of a ``hyperplane configuration''.

To explicitly describe geometric compactifications of the moduli space of $K3$ surfaces of a given degree  is a complicated and mostly open problem for many cases. We illustrate the case of degree equal to $2$, since generically those K3 surfaces are double covers of $\PP^2$ branched along a smooth plane sextic. Then one can be tempted to use the GIT quotient parametrizing such plane curves --- a technique used in the case of cubic surfaces. Indeed, GIT guarantees the existence of a semistable locus 
$\mathbb{P}H^0(\PP^2,\co(6))^{\text{ss}}$
such that the quotient
\begin{align*}
\overline{\M}^{\tiny{K3, \text{GIT}}}_{2}
:=
\mathbb{P}H^0(\PP^2,\co(6))^{\text{ss}}
/ \! \! /
\mathrm{SL}_3(\CC)
\end{align*}
is a compact projective variety of the correct dimension.  However, we should only think of 
$\overline{\M}^{\tiny{K3, \text{GIT}}}_2$ as a first approximation to our desired moduli space.  The reason is that $\M_2^{K3}$ contains an $18$-dimensional subfamily of \emph{unigonal} $K3$s whose map\footnote{
The degree-2 linear
system $|D|$ has a base curve $R$, and $D = 2E + R$. The free part of $|D|$ maps $X$ to the conic. }
to $\PP^2$ induced by the quasi-polarization is not surjective; instead, it produces an elliptic fibration of $X$ over a smooth conic $\{ (x_0x_1-x_2^2)=0 \}\subset\PP^2$.
The compact space $\overline{\M}^{\tiny{K3, \text{GIT}}}_2$ does not account for such $K3$s. Instead, it has an unique point $\omega$ parametrizing the sextic 
$\{ (x_0x_1 -x_2^2)^3 =0 \}$.  In 1980 Shah showed that the blow-up at $\omega$ yields a compact space
$ \Bl_{\omega}(\overline{\M}^{\tiny{K3, \text{GIT}}}_{2})$
containing $\M_2^{\tiny{K3}}$ as an open subset. The exceptional locus 
$\mathcal{E} \subset \Bl_{\omega}(\overline{\M}^{\tiny{K3, \text{GIT}}}_{2})$ 
is isomorphic to
\begin{align*}
\mathbb{P}\left(
 \operatorname{Sym}^{12}\left(\textbf{st}\right)
\times 
\operatorname{Sym}^{8}\left(\textbf{st} \right)
\right)
/ \! \! /_{\mathcal{O}(3,2)}
\mathrm{SL}_2(\CC)
\end{align*}
where $\textbf{st}$ is the standard $\mathrm{SL}_2(\CC)$ representation.  Geometrically, one interprets a general point in $\mathcal{E}$ as a pair of collections of 12 resp. 8 distinct points in $\PP^1$ up to an action of $\mathrm{SL}_2(\CC)$. To see the relation with unigonal $K3$ surfaces, note that the latter may be written in the form
$$
z^2 -y^3- yg_4(x_0,x_1,x_2) - g_6(x_0,x_1,x_2) = x_0x_1-x_2^2= 0.
$$
We recover collections of 12 resp. 8 points by considering the intersections $\{ g_4(x_0,x_1,x_2) = x_0x_1-x_2^2 = 0 \}$ and $\{ g_6(x_0,x_1,x_2) = x_0x_1-x_2^2 = 0 \}$.

Our visit to the case of points in the line in \S\ref{sec:pointsP1} suggests that there should be another type of compactification based on ``stable'' degenerations, a topic to be elaborated on in later sections.  Here we only remark that  the analogous compactification for degree-2 $K3$s is significantly more complicated than the GIT one,  and has only been completed in 2019 by Alexeev, Engel, and Thompson. 

The increasing complexity of the moduli phenomena makes some previously discussed questions, such as inversion of the period map, impracticable for arbitrary degree.  Therefore, a common strategy is to consider a smaller moduli problem for which one is able to give a more complete description. For example, one can consider either $K3$ surfaces whose Picard group is isometric to a particular lattice, or surfaces with a given finite automorphism group.  

From a Hodge-theoretic perspective, the \emph{vanishing} of periods on a sublattice $\LL\subset H_2(X,\ZZ)\cong H^{\oplus 3}\oplus E_8^{\oplus 2}$ of rank $\rho$ cuts out a subdomain $D_{\LL}\subset D$ of dimension $20-\rho$ (on which $\Gamma_{\LL}:=\text{stab}_{\Gamma_g}(\LL)$ acts, and which is always a ball or type-IV domain).  It follows from Lefschetz's (1,1)-theorem that its preimage $\phi^{-1}(D_{\LL}/\Gamma_{\LL})=:\M_{\LL}$ parametrizes \emph{$\LL$-polarized} $K3$s $(X,\ell)$, where $\ell\colon \LL\hookrightarrow \mathrm{Pic}(X)$ is an embedding of $\LL$ into the classes of line bundles (or divisors) on $X$.  For instance, for $\LL=H\oplus E_8\oplus E_7$ ($\rho=17$), Clingher and Doran have shown that the points 
$$[\alpha_2{:}\alpha_3{:}\alpha_5{:}\alpha_6]\in\M_{\LL}\cong \mathbb{WP}[2{:}3{:}5{:}6]\setminus \{\alpha_5{=}\alpha_6{=}0\}$$
parametrize $K3$s $X_{\ua}$ given by the minimal resolution of the quartic hypersurface in $\PP^3$ with affine equation
$$y^2z-4x^3z+3\alpha_2xz+\alpha_3z+\alpha_5xz^2-\tfrac{1}{2}(\alpha_6z^2+1)=0.$$
The period map $\phi$ restricts to an isomorphism 
$$
\M_{\LL}\overset{\cong}{\to}D_{\LL}/\Gamma_{\LL}\cong \hp_2/\mathrm{Sp}_4(\ZZ)
$$ 
which is inverted by the Siegel modular forms 
\footnote{
Here $\hp_2:=\{\btau\in M_2(\CC)\mid {}^t\btau=\btau,\,\mathrm{Im}(\btau)>0\}$ is Siegel's upper half-space, $\gamma=\left(\begin{smallmatrix}A&B\\C&D\end{smallmatrix}\right)\in \mathrm{Sp}_4(\ZZ)$ acts by $\gamma(\btau)=(A\btau+B)(C\btau+D)^{-1}$, and $f\in \mathcal{O}(\hp_2)$ belongs to $M_k(\hp_2,\mathrm{Sp}_4(\ZZ))$ iff $f(\gamma(\btau))=\{\det(C\btau+D)\}^k f(\btau)$.}
$\alpha_j\in M_{2j}(\hp_2,\mathrm{Sp}_4(\ZZ))$.

What about the \emph{blowing-up} of periods?  For a 1-parameter degeneration of $K3$ surfaces $\cx\to\Delta$, after a finite base-change ($t\mapsto t^m$) there are two possibilities: (a) two periods blow up like $\log(t)$; or (b) one period blows up like $\log(t)$, and another like $\log^2(t)$.  Geometrically, after birational modifications the singular fiber in (a) looks like two rational surfaces glued along an elliptic curve $E$ (or a more complicated variant of this), while in (b) we get a ``sphere-like'' configuration of rational surfaces glued along $\PP^1$'s (and the $\log^2(t)$ reflects the ``triple points'' in such a configuration).  In the compactification $\mathbb{WP}[2{:}3{:}5{:}6]\cong (\hp_2/\mathrm{Sp}_4(\ZZ))^*$, the point $[1{:}1{:}0{:}0]$ corresponds to (b), and the rest of $\{\alpha_5{=}\alpha_6{=}0\}$ to (a); the latter must be 1-dimensional to keep track of the modulus of $E$, after all.

As in the case of Picard curves, we can refine the ``limiting period map'' by computing additional periods that remain finite as $t\to 0$.  The idea is to look at the intersections of certain divisors on the rational surfaces with the rational or elliptic curves along which they are glued.  One then gets semiperiods by integrating (a) $\omega_E$ (on $E$) or (b) $\tfrac{dz}{z}$ (on $\PP^1\setminus \{0,\infty\}$) between the intersection points.  In the example, this contributes 1 modulus to (a) and 2 moduli to (b), making them into divisors and yielding a ``toroidal'' compactification $(\hp_2/\mathrm{Sp}_4(\ZZ))^{**}$.  (Of course, one does not get an extension of the period map \emph{to} this without first blowing up $\overline{\M}_{\LL}$).

\subsection{Cubic 4-folds}\label{sec:Cubic4folds}

This case is similar to the case of K3 surfaces.  A smooth cubic $X\subset \PP^5$ has (up to scale) a unique closed $(3,1)$-form, whose periods define a map $\phi$ from 
$\M:=\mathbb{P}\left(H^0(\PP^5,\mathcal{O}(3)\right)^{\text{sm}}/\mathrm{SL}_6(\CC)$
to the quotient $D/\Gamma$ of a 20-dimensional type-IV domain by the monodromy group.  By work of Laza, Looijenga and Voisin, $\phi$ is injective with image the complement of a hyperplane arrangement.  Cubics with action by a fixed finite group yield moduli subspaces, which behave analogously to the $\{\M_{\LL}\}$ in \S\ref{sec:K3surfaces}:  their images under $\phi$ turn out (by recent work of Laza-Pearlstein-Zhang and Yu-Zheng) to be hyperplane-arrangement complements in type-IV or ball quotients.

\subsection{Abelian varieties}

Given a full lattice $\Lambda :=\ZZ\langle \ula^{(1)},\ldots,\ula^{(2g)}\rangle\subset \CC^g$, the compact complex torus $T:=\CC^g/\Lambda$ has bases $\{dz_i\}_{i=1}^g$ of $\Omega^1(T)$ and $\{\gamma_j:=\overline{\uo.\ula^{(j)}}\}_{j=1}^{2g}$ of $H_1(T,\ZZ)$, hence $g\times 2g$ period \emph{matrix} $\Pi:=(\int_{\gamma_j}dz_i)=(\lambda_i^{(j)})$.  By Kodaira's embedding theorem, $T$ is algebraic\footnote{More concretely, $T$ is embedded in projective space by sections of powers of the theta line bundle.} (an \emph{abelian variety}) iff it possesses a closed, positive $(1,1)$-form $\Omega$ with $[\Omega]\in H^2(T,\ZZ)$. Equivalently, after changing both bases above, $\Pi$ takes the form $(\bm{\delta}\mid \btau)$, where $\bm{\delta}=\text{diag}(\delta_1,\ldots,\delta_g)\in M_g(\ZZ)$ (with $\delta_i\mid \delta_{i+1}$) and $\btau$ belongs to the Siegel upper half space $\hp_g:=\{\btau\in M_g(\CC)\mid {}^t\btau=\btau,\,\mathrm{Im}(\btau)>0\}$.

If all $\delta_i=1$, then $\Omega$ is a \emph{principal polarization} of $T$.  The set of principally polarized abelian varieties (PPAVs) is then parametrized by $\hp_g\ni \btau\mapsto A_{\btau}:=\CC^g/\Lambda_{\btau}$, where $\Lambda_{\btau}$ is the $\ZZ$-span of columns of $(\bm{1}_g\mid \btau)$. In order that each isomorphism class occur only once, we take the quotient $\A_g:=\hp_g/\mathrm{Sp}_{2g}(\ZZ)$, where $\gamma=\left(\begin{smallmatrix}\mathsf{A}&\mathsf{B}\\ \mathsf{C}&\mathsf{D}\end{smallmatrix}\right)$ acts by $\gamma(\btau)=(\mathsf{A}\btau+\mathsf{B})(\mathsf{C}\btau+\mathsf{D})^{-1}$.  This action on the ``period point'' is equivalent to a change of integral basis preserving $[\Omega]$ (viewed as a nondegenerate alternating form on $H_1(A_{\btau},\ZZ)$).  Siegel modular forms embed $\A_g$ into projective space, yielding our algebraic moduli space for PPAVs.  The compactification so obtained has a stratification of the form $\A_g^*=\A_g\amalg \A_{g-1}\amalg\cdots\amalg \A_0$.

Given a projective algebraic curve $C$ of genus $g$, the Jacobian $J(C)\cong \Omega^1(C)^{\vee}/H_1(C,\ZZ)$ is principally polarized by the intersection form on $H_1(C)$.  The resulting Torelli map $\Phi_g\colon \M_g\to \A_g$ (sending $C\mapsto J(C)$) is an injection, which is an immersion off the hyperelliptic locus.  Since the inequality (for $g>1$) $\dim_{\CC}\M_g=3g-3\leq {{g+1}\choose{2}}=\dim_{\CC}\A_g$ is strict for $g\geq 4$, most PPAVs are not Jacobians, and the period map for $\M_g$ does \emph{not} have open image.

\section{Aspects of the theory}

Our tour has reached the base of the funicular, on which we now ascend to two different vantage points for a brief theoretical overview.  On the algebro-geometric side, we refer the reader to \cite{kollar2010moduli}, \cite{xu2020k} and references therein for the technical details; a good place to start for curves are \cite{caporaso}, \cite{hassett}
and for higher dimensional cases  are \cite{alexeev2013moduli} and \cite{Kov05}.
For background on period maps in Hodge theory, we recommend \cite{CMSP}.

\subsection{What is a moduli space?}
In moduli theory, we want more than a space parametrizing objects, such as smooth elliptic curves up to isomorphism.  We want to understand all possible families of the objects.  For that purpose, we need to reformulate our problem.

Let $\Omega$ be a “reasonable” class of objects --- for example, the stable $n$-pointed curves of genus $0$ described in \S\ref{sec:pointsP1}. The corresponding \emph{moduli functor}\footnote{For the newcomer: ``flatness" ensures reasonable behavior, such as constancy of the Hilbert polynomial in a projective family; if $B$ is regular and $\mathcal{X}$ Cohen-Macaulay, it is simply equivalent to equidimensionality of the fibers.}
\begin{align*}
\mathcal{M}(B) &= 
\{
\text{flat families  $\mathcal{X} \to B$ 
 with fibers
 $X_b\in \Omega$,} 
\\
&
\hspace{1cm}
\text{
modulo isomorphisms  over $B$} 
\}    
\end{align*}
maps schemes $B$ to sets $\mathcal{M}(B)$.  For any scheme $\M$, we can define another functor $h_{\M}$ from schemes to sets by $h_{\M}(B) :=\mathrm{Mor}( B, \M)$.  A moduli functor $\mathcal M$ is \emph{represented} by a scheme $\M$ if there is a natural isomorphism from 
$\mathcal{M}$ to  $h_{\M}$.  In that event, $\M$ is called a \emph{fine moduli space} for $\M$, and constructing a family of objects over a base $B$ is equivalent to defining a morphism $B \to \M$.  

Looking back at the examples in \S\ref{sec:ExampleModuliTheory},  $\overline{\M}_{0,n}$ 
represents the moduli functor associated to stable $n$-pointed curves of genus $0$.  Another case in point is the Hilbert scheme, which represents the functor
\begin{align*}
\mathcal{M}(B) &= 
\{
\text{flat families  $\mathcal{X} \to B$ 
 such that 
 $X_b \subset \PP^r$ 
 } 
\\
&
\hspace{0.5cm}
\text{
has Hilbert polynomial equal to $p(m)$} 
\}.
\end{align*}
Unfortunately, most moduli functors are not represented by a scheme. Even the moduli functor associated to isomorphism classes of smooth elliptic curves fails to be represented by $\M_{1,1}$.

One option for moving forward is to weaken our expectations.  A scheme $\M$ is a 
\emph{coarse moduli space} 
for $\mathcal M$ if $\M$ is the scheme best approximating the moduli functor, and the geometric points of $\M$ are in bijection with the equivalence-classes of objects ``parametrized'' by $\mathcal{M}$. 
More precisely, there is a natural transformation $\mathcal{M}\longrightarrow h_{\M}$, for which $\M$ is universal and $\mathcal{M}(\Spec (K))\longrightarrow h_{\M}(\Spec (K))$ an isomorphism of sets in $K=\bar{K}$.  But for $B$ more general than a geometric point, $\mathcal{M}(B)\longrightarrow \mathrm{Mor}(B,\M)$ may be neither injective nor surjective:  distinct families over $B$ may produce the same ``classifying map'' to $\M$, and not every map to $\M$ (even $\mathrm{id}_{\M}$) need give rise to a family. On the positive side, if a moduli functor has a coarse moduli space, then the latter is unique (up to canonical isomorphism).  For example, the GIT quotient $\M_{1,1}$ is \emph{the} coarse moduli space for smooth elliptic curves up to isomorphism, though it neither carries a universal family nor distinguishes (say) isotrivial families from trivial ones.

 What if, in the absence of a fine moduli space, we still want to keep track of the families?  
 In that case, we need new geometric objects known as \emph{stacks}. Introduced by Deligne, Mumford, and Artin in the 1970s, they are (loosely speaking) enrichments of schemes obtained by attaching an automorphism group to every point. 
In our particular context, the stack of objects 
$\Omega$ is a category whose objects are families $\{\mathcal{X} \to T\}$ of our ``reasonable" objects. A morphism in this category 
$\{\mathcal{X} \to T\} \mapsto \{ \mathcal{Y} \to B \}$ 
is a pair of maps $f: \mathcal{X} \to \mathcal{Y}$ and 
$g: T \to B$ such that the diagram
\begin{align*}
    \xymatrix{
    \mathcal{X} \ar[r]^f \ar[d] &  \mathcal{Y} \ar[d]
    \\
    T \ar[r]^g                  & B
    }
\end{align*}
commutes and $\mathcal{X}$ is isomorphic to the pullback of $\mathcal Y$ via the map $T \mapsto B$.  

These two roads, stacks and coarse moduli spaces, converge to yield a rich picture.  For a well-behaved moduli problem (for example, when automorphisms of all parametrized objects are finite and the limits are uniquely determined, i.e. the stack is separated), Keel and Mori showed in 1997 that there is a coarse moduli space $\M$ associated with the stack describing our moduli functor $\mathcal{M}$.  
Recent results by Alper, Halpern-Leistner, and Heinloth have generalized this result to a larger class of stacks, whose parametrized objects can have positive-dimensional automorphism groups.  A caveat here is that the output $\M$ is an algebraic space rather than a scheme.  

The previous discussion leaves us with a delicate question. What are the ``reasonable" classes of objects for compactifying a given moduli space?  Or, more colloquially:  \emph{what do we add in the boundary?}  The answer depends on the types of varieties we are considering.  We now visit two important cases.

\subsection{\normalsize Algebro-geometric compactifications} 
\label{sec:geocompact}

\subsubsection{General type case}
Let $\M$ be a moduli space of pairs $(X,D)$ where $X$ is smooth projective, $D$ is a normal crossing \emph{effective} divisor, and $K_X +D$ is ample.  Examples include (non-Eckardt) smooth cubic surfaces together with their $27$ lines, or complex curves of genus $g \geq 2$ (with $D=0$).  In terms of the minimal model program, the ``correct'' objects  for compactifying such $\M$ are the (more general) \emph{stable pairs}.  These consist of a projective variety $X$ and an $\RR$-divisor $D=\sum b_i D_i$ on it, such that $K_X+D$ is $\RR$-Cartier and ample, and the pair has semi-log-canonical (slc) singularities, see \cite[Def 1.3.1]{alexeev2013moduli}.   

One reason to add a divisor to $X$, is that $K_X$ may not be ample; in that case, without the ``boost'' from $D$, the canonical model wouldn't recover our original family, let alone compactify it.  On the other hand, the choice of $D$ may add dimensions to our moduli space; or there may be various choices of $D$ that lead to distinct compactifications of the same $\M$.

We first illustrate the concepts surrounding slc singularities when $D=0$ and $X$ is singular but normal. Let $f\colon Y \to X$ be a resolution of singularities such that the exceptional set $\cup \mathcal{E}_j$ is a normal crossing divisor. Assume that the canonical divisor $K_X$ is  $\QQ$-Cartier; then we have the numerical equivalence
\begin{equation}\label{eq-disc}
K_Y\sim_{\QQ} f^{*}(K_X) + \textstyle{\sum_j} a_j \mathcal{E}_j
\end{equation}
of $\QQ$-Weil divisors.  The coefficients $a_j$ are called the \emph{discrepancies}, and they measure how bad our singularities are.  The singularities of $X$ are called \emph{terminal} (the best case) if $a_j >0$, \emph{canonical} (still very good) if $a_j \geq 0$ and \emph{log-canonical} (tolerable) if $a_j \geq -1$ for all $j$.  In dimension 2, the canonical hypersurface singularities are precisely the ADE ones, while the log-canonical ones include all quotient singularities $\frac{1}{m}(1,r):=\mathbb{C}^2/\mu_m$ with $\gcd(m,r)=1$; $\mu_m$ is the group of $m^{\text{th}}$-roots of unity acting via $\rho(x,y) \mapsto (\rho x, \rho^r x)$ with $\rho^m=1$.  

The defintion of discrepancies extends to pairs with $D\neq 0$ by replacing $K_X$ in \eqref{eq-disc} by $K_X+D$.  We say a pair is \emph{klt} (resp. log-canonical) if all its discrepancies are strictly larger than $-1$ (resp. $ \geq1 $). The \emph{slc} singularities are those that become ``log-canonical after normalization."  They include (for $D=0$) the pinch point $\{ x^2 + y^2z =0 \}$ and degenerate cusps $\{ y^2\left( z^2 + y^{q-2}\right) = 0 \}$ (with $q \geq 3$).

Even if we know that stable pairs are the right objects for compactifying our moduli, defining \emph{families} of them is a tricky business that leads to two possible functors, due to 
Koll\'{a}r
and Viehweg.  Moreover, if $D \neq 0$, the ``boundary fibers'' may exhibit phenomena like non-reduced points, see \cite[Sec 1.1]{KZ}.  So here we shall limit ourselves to an oversimplified\footnote{see \cite[\S1.5]{alexeev2013moduli} and \cite[Sec 5]{KZ} for a more precise statement}
definition of the moduli functor that produces the so-called \emph{slc} or \emph{KSBA compactification} 
(after Koll\'ar, Shepherd-Barron, and Alexeev)  when $D=0$:
\begin{align*}
\mathcal{M}^{\tiny\text{slc}}(B) &= 
\{\text{flat proper morphisms  $\mathcal{X} \to B$ whose}\\
& \hspace{0.5cm}\text{ fibers $X_s$ are $n$-dimensional stable } \\
& \hspace{0.5cm}\text{  varieties with fixed  $K_{X_s}^{\dim(X)} \in \QQ_{>0}$,} \\ 
& \hspace{0.5cm}\text{ satisfying Koll\'{a}r's condition
}
\}  
\end{align*} 
This functor\footnote{
Here ``Koll\'ar's condition" requires that for any $m \in \mathbb{Z}$ the reflexive power $\omega_{\mathcal{X}/S}^{[m]}$ commutes with arbitrary base change.} is in fact coarsely represented by a projective scheme. This last claim is also true for the generalizations involving non-zero divisors $D$, see \cite[Corollary 6.3]{KZ} and \cite[Sec 1.6]{alexeev2013moduli}.

We now return to one of the main landmarks in algebraic geometry, of which the reader caught a glimpse in \S\ref{sec:pointsP1}.  Fix $n$ real numbers $0 < b_i \leq 1$. A  \emph{weighted stable curve} for the weight $\mathbf{b}= (b_1, \ldots, b_n)$ is a \emph{marked} pair (in the sense of \cite{kollar2010moduli}) $(C,(D_1,\ldots,D_n))$ comprising a reduced connected projective curve $C$ together with a collection of $n$ points $D_i\in C$, such that the divisor $K_C + D:=K_C+\sum_i b_i D_i$ is ample. The slc condition translates into: 
\begin{itemize}
    \item $C$ is either smooth or has at worst singularities locally analytically isomorphic to $\{ xy=0 \}$, and
    \item  the points $D_i$ may coincide with each other, but they should be different from the nodes, and
    \item  the sums of multiplicities of coincident points must not exceed $1$:  $\sum_{i\colon D_i=p} b_i\leq 1$ ($\forall p\in C$).
\end{itemize}
It  follows from works of Deligne, Mumford, Knudsen, and Hassett that  for any $n$, $\mathbf{b}$, and $g \geq 0$ the moduli stack  $\mathcal{M}_{g, \mathbf{b}}$ of weighted stable curves of arithmetic genus $g$ is a smooth Deligne-Mumford stack with a projective coarse moduli space $\overline{\M}_{g, \mathbf{b}}$.

\begin{figure}[h!]
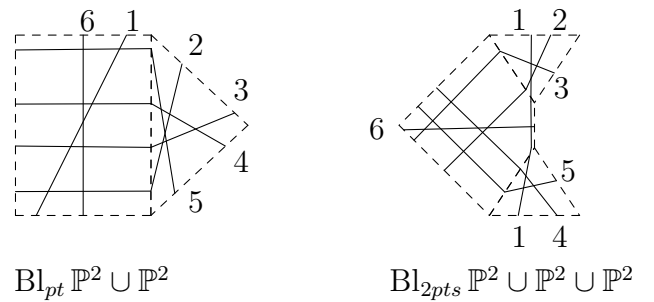

\centering
\tikzpicture[line cap=round,line join=round, 
scale=0.6]
\begin{scope}[shift ={(1,0)}]
\draw [dashed] (2,6)-- (2,2);
\draw [dashed] (2,2)-- (5,2);
\draw [dashed] (5,2)-- (5.01,6);
\draw [dashed] (5.01,6)-- (2,6);
\draw [dashed] (5.01,6)-- (7.15,4);
\draw [dashed] (7.15,4)-- (5,2);
\draw [dashed] (5,2)-- (5.01,6);

\draw [] (3.5,6)-- (3.5,2);
\draw []
(2.46,2)-- (4.46,6);
\draw []
(2,5.67)-- (5.01,5.7);
\draw []
(2,4.47)-- (5.01,4.48);
\draw []
(2,3.54)-- (5,3.51);
\draw []
(2,2.53)-- (5,2.54);
\draw []
(5,2.54)-- (5.69,5.36);
\draw []
(5.01,5.7)-- (5.53,2.49);
\draw []
(5,3.51)-- (6.85,4.27);
\draw []
(5.01,4.48)-- (6.65,3.54);

\draw  (3.60,6.32) node {$6$};
\draw (3.7,0.5) node {$\Bl_{pt}\PP^2 \cup \PP^2$ };

\draw  (4.6,6.32) node {$1$};
\draw (6.0,5.79) node {$2$};
\draw (7.00,4.79) node {$3$};
\draw (7.00,3.19) node {$4$};
\draw (6.0,2.3) node {$5$};
\end{scope}

\begin{scope}[shift ={(0.5,0)}]
\draw [dashed] (13,6)-- (11,4);
\draw [dashed] (11,4)-- (13,2);
\draw [dashed] (13,2)-- (14,3.5);
\draw [dashed] (14,3.5)-- (14,4.5);
\draw [dashed] (14,4.5)-- (13,6);
\draw [] (11.1,3.9)-- (14,3.97);
\draw []
(11.33,3.67)-- (13.24,5.64);
\draw []
(12,3)-- (13.81,4.79);
\draw []
(14.36,6)-- (13.81,4.79);
\draw []
(13.24,5.64)-- (14.44,5.16);
\draw []
(13.92,4.52)-- (13.93,3.51);
\draw []
(13.64,2)-- (13.93,3.51);
\draw []
(13.92,4.52)-- (13.93,6);
\draw []
(11.42,4.42)-- (13.34,2.52);
\draw []
(13.34,2.52)-- (14.48,2.77);
\draw []
(11.83,4.83)-- (13.69,3.03);
\draw []
(13.69,3.03)-- (14.5,2);

\draw [dashed] (13,6)-- (15,6);
\draw [dashed] (15,6)-- (14,4.5);
\draw [dashed] (14,4.5)-- (13,6);
\draw [dashed] (14,3.5)-- (15,2);
\draw [dashed] (15,2)-- (13,2);
\draw [dashed] (13,2)-- (14,3.5);

\draw (10.51,3.92) node {$6$};
\draw (13.5, 0.5) node 
{$\Bl_{2pts}\mathbb{P}^2 \cup \PP^2 \cup \PP^2$};

\draw (13.66,6.33) node {$1$};
\draw (14.56,6.33) node {$2$};
\draw (14.60,4.83) node {$3$};
\draw (13.66,1.53) node {$1$};
\draw (14.56,1.53) node {$4$};
\draw (14.75,3.0) node {$5$};
\end{scope}
\endtikzpicture
\caption{Examples of stable pairs}
\label{ref:PairsDim2}
\end{figure}

The explicit description of higher dimensional pairs parametrized by compact moduli spaces is a busy industry nowdays. 
A quite combinatorially rich example is the moduli of hyperplane arrangements in projective space \cite{alexeev2013moduli}.
We illustrate in Fig. \ref{ref:PairsDim2} two stable pairs (from [op. cit.]) found in the compactification of the moduli space of six labelled lines in $\PP^2$.


\subsubsection{Fano case} 
The above theory does not work when $K_X +D$ is not ample, so a new perspective is necessary.  We will discuss the case where $X$ is Fano (i.e. $-K_X$ is ample) and $D=0$.  Many of the ideas can be extended to the case when $-K_X-D$ is ample, and we encourage the reader to consult \cite{xu2020k} for the  details. 

Here, the concept of $K$-(semi)stability is central. Introduced in 1997 by Tian, it became a central conjecture --- now theorem --- that the existence of a KE metric on a smooth Fano variety $X$ is equivalent to satisfying a K-stability condition. 
Key ingredients are test configurations for $(X,L)$, i.e. certain one-parameter degenerations $\mathcal{X} \to \mathbb{C}$ with fibers $X_t \cong X$ ($\forall t\in \CC^*$) and a $\mathbb{C}^*$-action, and the Donaldson-Futaki (DF) invariant 
of such test configurations. A Fano variety $X$ is \emph{$K$-semistable} iff the DF invariant is non-negative for all test configurations for $(X,L)$. 

It is hard to check $K$-(semi)stability from its definition. However, there is a local-to-global interplay that restricts the geometry of $K$-semistable varieties.  By work of Odaka, a reasonable\footnote{$X$ is equidimensional, reduced, $S_2$, the codimension one points of $X$ are Gorenstein, and $K_X$ is $\QQ$-Cartier} $K$-semistable Fano variety has at worst klt singularities. 
Moreover, Liu showed in 2018 that for an $n$-dimensional $K$-semistable Fano variety $X$, a local invariant at $p \in X$, denoted by $\widehat{vol}(p,X)$,\footnote{
Let $(X, p)$ = $(Spec(R), \mathbf{m})$, where $R$ is a local ring essentially of finite type and $\mathbf{m}$ is the maximal ideal; then
$
\widehat{vol}\left( p, X\right)
=
\inf \left\{
lct(\mathbf{a})^n\operatorname{mult}(\mathbf{a})
\; | \; 
\mathbf{a} \text{ is } \mathbf{m}-primary
\right\}
$
}
is bounded by global 
invariants of the variety: 
$
\widehat{vol}\left( p, X\right)
\left( \frac{n+1}{n} \right)^n
\geq (-K_X)^n. 
$
These relations play a central role in finding  singular $K$-semistable varieties. 

The construction of the moduli space of $K$-semistable varieties is ongoing work of many people including Alper, Blum, Halpern-Leistner, Li, Liu, Wang, Xu, Zhuang, among others. The moduli functor $\mathcal{M}^{\KE}$ sends the scheme $B$ to 
\begin{align*}
\mathcal{M}^{\KE}(B) &= 
\{\text{flat proper morphisms $\mathcal{X} \to B$ whose} \\
&\hspace{0.5cm}\text{fibers are $n$-dimensional $K$-semistable} \\
&\hspace{0.5cm}\text{varieties with volume $\mathrm{V}_0$, satisfying}\\
&\hspace{0.5cm}\text{Koll\'ar's condition}\}.  
\end{align*}
This functor is represented by an Artin stack of finite type and admits a projective coarse moduli space $\overline{\M}^{\KE}$.  
We remark that it is possible to extend the above concepts to pairs. 
Particular cases are described by work of the first author with Martinez-Garcia and Spotti; while a more general setting is given by 
Ascher, Liu and DeVleming.  But this is another story and shall be told another time.

\subsubsection{A brief stop at GIT}
If our moduli space parametrizes, say, cubic surfaces up to isomorphism, shouldn't the limit also be a cubic surface? KSBA or $K$-stability techniques do not guarantee or even prioritize such a condition. However,  Geometric Invariant Theory (GIT) does.

The main idea is that, given a moduli space $\M$ for smooth varieties or pairs, we can select an embedding that realizes such objects inside a projective space $\PP^{N-1}$. We can then parametrize those embedded objects by an open subset $\mathcal{U}$ in a Hilbert scheme (or proxy).  
Under appropriate circumstances, 
there exists a line bundle $L$ with an $\SL_N(\mathbb{C})$-linearization, and a larger open $\mathcal{U}^{\text{ss}}(L)\,(\supset \mathcal{U})$ called the \emph{semistable locus} (depending on $L$), such that the categorical quotient $\mathcal{U}^{ss}(L)/ \! \!/ \SL_{N}(\mathbb{C})$ is a projective variety compactifying our moduli space $\M$.  Moreover, there is an intermediate \emph{stable locus} $\mathcal{U}^{\text{s}}(L)$ (also depending on $L$) for which the geometric quotient  $\mathcal{U}^s(L)/\SL_{N}(\mathbb{C})$ is well-behaved, and
$$\M\subseteq \mathcal{U}^{\text{s}}(L)/\SL_N(\CC)\subseteq \mathcal{U}^{\text{ss}}(L)/\!\!/\SL_{N}(\CC).$$
We recover our initial moduli problem because $\M \cong \mathcal{U}/\SL_{N}(\mathbb{C})$. 

The choice of the embedding and of the line bundle $L$ mean that we have many possible GIT compactifications for $\M$. In particular, Gieseker showed that if we use the pluricanonical embedding with $m \geq 5$  for complex curves of genus $g \geq 2$, 
then we can construct a GIT quotient that recovers $\overline{\M}_{g}$. Other GIT compactifications of $\M_g$ and $\M_{g,\mathbf{b}}$ have been the subject of work by Hassett, Hyeon, Fedorchuk, Lee, Li, Schubert, Swinarski,  Smyth, Jensen,  Moon, and others.  
A different example where multiple GIT quotients come into play, this time involving cubic surfaces, was described by the first author with Martinez-Garcia \cite{GMG}.

\subsubsection{Remarks on other cases}

We can use some of the above techniques for compactifying moduli of varieties for which neither $K_X$ nor $-K_X$ are ample.  For instance, if $X$ is a Calabi-Yau variety, we can select an ample divisor $D$, so that the pair $(X,\epsilon D)$ is stable.  Koll\'ar and Xu have recently determined that the irreducible components of the resulting moduli space are projective.  Other compactifications may be obtained by using either GIT or the GKZ secondary fan when $X$ is a complete intersection in a toric variety.

\subsection{\normalsize Hodge-theoretic compactifications}\label{sec:HodgeTheory}

\subsubsection{Period maps}

A weight-$n$ \emph{Hodge structure} on a $\QQ$-vector space $V$ is a decomposition $V_{\CC}=\oplus_{p+q=n}V^{p,q}$ with $\overline{V^{q,p}}=V^{p,q}$.  By the Hodge theorem, the classes $H^{p,q}(X)\subset H^n(X,\CC)$ of closed $C^{\infty}$ forms of type $(p,q)$ yield a HS for any compact K\"ahler manifold $X$.  Equivalent data are given by the decreasing filtration $F^{\bullet}$ defined by $F^k V=\oplus_{p\geq k} V^{p,n-p}$ or homomorphism 
$$
\vf\colon  S^1=\{z\in \CC\bigl\vert |z|=1\}\to \mathrm{SL}(V)
$$ 
defined by $\vf(z)|_{V^{p,q}}=z^{p-q}.\mathrm{Id}_{V^{p,q}}$. The \emph{Hodge numbers} of $V$ are $\uh=(h^{p,q}):=(\dim_\CC V^{p,q})$, and its \emph{level} is $\max\{|p-p'|\bigl\vert h^{p,n-p}\neq 0\neq h^{p',n-p'}\}$.

A \emph{polarization} of the HS $(V,F^{\bullet})$ is a nondegenerate, $(-1)^n$-symmetric bilinear form $Q\colon V\otimes V\to \QQ$ satisfying the Hodge-Riemann bilinear relations
\begin{subequations}
\begin{align}
Q(F^p,F^{p'})=0\;&\text{ if $p+p'>n$}\label{I4a} \\
\sqrt{-1}^{p-q}Q(v,\bar{v})>0\;&\text{ if $v\in V^{p,q}\setminus\{0\}$}.\label{I4b}
\end{align}
\end{subequations}
If the manifold $X$ is projective, then $H^n(X)$ admits a polarization by the hard Lefschetz theorem.

The classifying spaces $\D_{(\uh,Q)}$ for HSs on $V$ with given $\uh$ and polarized by $Q$ are called \emph{period domains}.  Writing $f^p:=\sum_{p'\geq p}h^{p',n-p'}$, \eqref{I4a} defines a projective subvariety $$\check{\D}_{(\uh,Q)}\subset \textstyle{\prod_p} \mathrm{Grass}(f^p,V_{\CC})$$ called the \emph{compact dual}, inside which \eqref{I4b} describes $\D_{(\uh,Q)}$ as an analytic open subset.  Writing $\mathcal{G}:=\mathrm{Aut}(V,Q)$ (an orthogonal resp. symplectic $\QQ$-algebraic group for $n$ even resp. odd), $\mathcal{G}(\RR)$ acts transtitively on $\D_{(\uh,Q)}$ and $\mathcal{G}(\CC)$ on $\check{\D}_{(\uh,Q)}$.  There is a plethora of (\emph{Hodge} or \emph{Mumford-Tate}) subdomains $$D\subset \D_{(\uh,Q)},$$ constructed as follows:  the \emph{Hodge group} of a PHS $(V,F^{\bullet},Q)$ is the smallest $\QQ$-algebraic subgroup $G\leq \mathcal{G}$ with $\vf(S^1)\leq G(\RR)$, and we set $D:=G(\RR)^+.F^{\bullet}$ (and $\check{D}:=G(\CC).F^{\bullet}$). For instance, these may parametrize HSs stabilized by a fixed cyclic automorphism of $V$, as we saw in several examples in \S\ref{sec:ExampleModuliTheory} (with $D$ a complex ball and $G(\RR)$ indefinite unitary).

Now let $\{X_{\sm}\}_{\sm\in \M}$ be a smooth algebraic family of smooth projective varieties. We may identify $H^n(X_{\sm},\QQ)$ with a fixed $V$ up to the action of monodromy, and consider the variation of the structures on $V_{\CC}$ induced by the Hodge theorem.  Write $D$ for the smallest Hodge subdomain containing all these structures, and $\Gamma\leq G(\QQ)$ for the monodromy group.  The key observation, due to Griffiths, is that while $V_{\sm}^{p,q}$ varies \emph{non}holomorphically, $F^p_{\sm}V_{\CC}$ varies holomorphically and satisfies the \emph{infinitesimal period relation} (a.k.a. \emph{Griffiths transversality})
\begin{equation}\label{I5}
dF^p_{\sm}\subset \Omega^1_{\M}\otimes F^{p-1}_{\sm}.
\end{equation}
Abstractifying this yields the notion of a (polarized) \emph{variation} of HS $\mathcal{V}\to\M$, consisting of a $\QQ$-local system $\mathbb{V}$ and a filtration of $\mathbb{V}\otimes \mathcal{O}_{\M}$ by (sections of) holomorphic subbundles, such that the fibers are PHSs, and the flat connection annihilating $\mathbb{V}$ satisfies \eqref{I5}.

Any PVHS yields a (locally liftable, holomorphic) \emph{period map}
\begin{equation}\label{I6}
\phi\colon \M\to D/\Gamma,
\end{equation}
where $\phi(\M)$ is an integral manifold of the $G(\RR)$-invariant horizontal distribution $\W\subset TD$ equivalent to \eqref{I5}.  We want to use $\phi$, together with a (possibly partial) compactification of $D/\Gamma$, to obtain an algebraic compactification of $\M$ with (mixed) Hodge-theoretic modular meaning.  This seems most plausible if \emph{(i)} $\phi$ is injective (i.e. a \emph{global Torelli theorem} holds), \emph{(ii)} $D/\Gamma$ is algebraic, and \emph{(iii)} $\phi$ is onto a Zariski open subset --- as in the examples of \S\ref{sec:ExampleModuliTheory}.  The bad news is that, if $\W$ is not $TD$, then it obstructs \emph{(iii)}, and $D/\Gamma$ is usually non-algebraic (Griffiths-Robles-Toledo).  Though the situation is better than that makes it sound, we will start by discussing compactifications in the ``classical'' case where $\W=TD$, which is to say, where the family of HSs parametrized by $D/\Gamma$ satisfies \eqref{I5} (yielding a \emph{tautological VHS}).

\subsubsection{The classical case}

Given a HS $\vf\in D$, $\vf(z)$ acts on the Lie algebra $\fg_{\CC}\subset \mathrm{End}(V_{\CC})$ of $G(\CC)$ by conjugation, inducing a Hodge decomposition $\fg_{\CC}=\oplus_p \fg_{\vf}^{p,-p}=:\oplus_p \fg^p_{\vf}$ of weight $0$.  Since $T_{\vf}D$ is identified with $\oplus_{p<0} \fg_{\vf}^p$ and $\W_{\vf}$ with $\fg_{\vf}^{-1}$, we see that $\W=TD$ forces $\fg_{\vf}^p=\{0\}$ for $p\neq -1,0,1$.  So the ``symmetry'' $\vf(\sqrt{-1})\colon D\to D$ fixing $\vf$ has differential $-\mathrm{Id}_{T_{\vf}D}$ at $\vf$, and $D\cong G(\RR)^+/G^0(\RR)\cong G^{\text{ad}}(\RR)^+/K$ is a \emph{Hermtian symmetric domain} (of noncompact type). Conversely, the irreducible HSDs are classified by special nodes on connected Dynkin diagrams, which leads to the list
\begin{itemize}[leftmargin=0.5cm]
\setlength\itemsep{0em}
\item $\mathsf{A}_n$:  $D\cong \mathrm{SU}(p,q)/\mathrm{S}(\mathrm{U}(p){\times}\mathrm{U}(q))$ 

$(p{+}q=n{+}1;\; p,q>0)$
\item $\mathsf{B}_n$: $D\cong \mathrm{SO}(2,2n{-}1)^+/(\mathrm{SO}(2){\times}\mathrm{SO}(2n{-}1))$
\item $\mathsf{C}_n$: $D\cong \mathrm{Sp}_{2n}(\RR)/\mathrm{U}(n)$
\item $\mathsf{D}_n$: $D\cong \mathrm{SO}(2,2n{-}2)^+/(\mathrm{SO}(2){\times}\mathrm{SO}(2n{-}2))$ 

or $\mathrm{SO}^*(2n)/\mathrm{U}(n)$
\item $\mathsf{E}_6,\mathsf{E}_7$: exceptional HSDs.
\end{itemize}
They are realized as Hodge domains (in some connected component $\D_{\uh}^+$ of a period domain) by taking a representation $V$ of $G$, and considering the tautological VHS (of some weight $n$) induced by the decomposition $\fg_{\CC}=\fg_{\vf}^{-1}\oplus\fg_{\vf}^0\oplus\fg_{\vf}^1$ (and polarized by the Killing form) as $\vf$ varies over $D$.  In particular, taking $V$ to be the standard representation in case $\mathsf{B}_{\frac{k+1}{2}}$/$\mathsf{D}_{\frac{k}{2}+1}$ resp. $\mathsf{C}_g$ yields $D=\D^+_{(1,k,1)}$ resp. $\D_{(g,g)}(\cong \hp_g)$; and these are the only full period domains with $\W=T\D$.  Of course, for Hodge numbers different from these (i.e. $(a,b,a)$ with $a\geq 2$, or level $\geq 3$), we can still have ``classical'' period map targets; they just have to be \emph{proper} Hodge subdomains.

Now assume that $\Gamma \leq G(\QQ)$ stabilizes a lattice $V_{\ZZ}\subset V$ (as monodromy groups do), and that furthermore $\Gamma$ is of finite index in $G(\QQ)\cap \mathrm{Aut}(V_{\ZZ})$ (i.e. $\Gamma$ is \emph{arithmetic}).  The Baily-Borel theorem then guarantees the existence of enough $\Gamma$-invariant sections of powers of the canonical bundle\footnote{One can also use the \emph{Griffiths line bundle} $L_{D/\Gamma}:=\otimes_p (\det V^{p,n-p})^{\otimes p}$, which is always a rational tensor power of $K_{D/\Gamma}$ in the classical case.} $K_D$ (i.e. modular forms) to present $D/\Gamma$ as a quasi-projective variety, which is smooth for $\Gamma$ torsion-free (and called a \emph{locally symmetric variety}).  The resulting compactification $$(D/\Gamma)^*:=\mathrm{Proj}(\oplus_k H^0(D/\Gamma,K_{D/\Gamma}^{\otimes k}))$$ is (as a set) a disjoint union of $D/\Gamma$ and finitely many \emph{Baily-Borel (B-B) boundary strata} $B^*/\Gamma_{B^*}$, where $B^*$ runs over $\Gamma$-equivalence classes of rational holomorphic path components in $\partial D\subset \check{D}$.

These strata have a Hodge-theoretic characterization.  Given a period map \eqref{I6} with locally symmetric target, and $\overline{\M}\supset \M$ any good compactification,\footnote{This means that $\overline{\M}$ is smooth and proper, and $\overline{\M}\setminus \M$ is a simple normal crossing divisor.} \emph{Borel's extension theorem} yields a holomorphic map $\phi^*\colon \overline{\M}\to (D/\Gamma)^*$ restricting to $\phi$.  (By GAGA, $\phi^*$ hence $\phi$ is algebraic; and it follows that the above algebraic structure on $D/\Gamma$ is unique.)  The key point is then that $\phi^*|_{\overline{\M}\setminus \M}$ \emph{simply records the limit of the flag} $F^{\bullet}\subset V$, up to $\Gamma$-equivalence.  This has the effect (seen in the examples of \S\ref{sec:ExampleModuliTheory}) of killing off all finite limits of periods of a form $\omega$ if one of its periods blows up, since these periods are projective coordinates of a line $\CC[\omega]\subseteq F^n\subset V$.

Hodge-theorists express limits of periods in terms of \emph{limiting mixed Hodge structures (LMHS)}. Passing to a finite cover, we may assume $\Gamma$ is neat.\footnote{i.e., its eigenvalues generate a torsion-free subgroup of $\CC^*$.}  Let $z_1,\ldots,z_d$ be local coordinates at $x\in \overline{\M}\setminus\M\overset{\scriptscriptstyle{\text{loc}}}{=}\{z_1\cdots z_k=0\}$.  Then the corresponding local monodromy generators $T_1,\ldots,T_k\in \Gamma$ commute and are unipotent, so that the $N_i=\log(T_i)$ generate a rational nilpotent cone $\sigma=\sigma_x\subset \fg_{\RR}$.  There is a unique increasing filtration $W_{\bullet}$ of $V$ satisfying $NW_{\bullet}\subset W_{\bullet-2}$ and $N^k\colon \mathrm{Gr}^W_{n+k}V\overset{\cong}{\to}\mathrm{Gr}^W_{n-k}V$ for any $N\in \mathrm{int}(\sigma)$.  Setting
$$F_{\lim}^{\bullet}:=\lim_{\uz\to \uo} e^{-\sum_i \frac{\log(z_i)}{2\pi\sqrt{-1}}N_i}F_{\uz}^{\bullet},$$ the LMHS at $x$ is the triple $(V,W_{\bullet},F_{\lim}^{\bullet})$, which we record modulo $e^{\CC\sigma}$ (to eliminate the dependence on the $\{z_i\}$) and $\Gamma_{\sigma}:=\mathrm{stab}(\sigma)\cap \Gamma$ in a \emph{Hodge-theoretic boundary component} $B(\sigma)/\Gamma_{\sigma}$, see \cite{KP}.  Passing to the $\mathrm{Gr}^W_{\bullet}$ of the LMHS maps this down to a B-B stratum, forgetting ``extension classes'' (like the semi-periods in the examples).

Given a suitable $\Gamma$-compatible collection of such cones (closed under taking intersections and faces) called a \emph{projective fan} $\Sigma$ in $\fg_{\RR}$, $D/\Gamma$ and the $\{B(\sigma)/\Gamma_{\sigma}\}_{\sigma\in \Sigma}$ glue together into a smooth projective \emph{toroidal compactification} $(\overline{D/\Gamma})_{\Sigma}$, which is a resolution of singularities of $(D/\Gamma)^*$ \cite{AMRT}.  Such a fan always exists.\footnote{For ball quotients, the cones are rays and it is unique.  On the other hand, for quotients of $\hp_g$, the reduction theory of quadratic forms provides several options (cf. \S\ref{sec:abelian}).} In general, $\phi^*$ factors through $(\overline{D/\Gamma})_{\Sigma}$ only after birationally modifying $\overline{\M}$ along $\overline{\M}\setminus\M$.

Finally, if $D$ is a ball or type-IV domain, and $\mathcal{H}$ a $\Gamma$-invariant hyperplane configuration, there is an important \emph{semi}toroidal compactification $({(D\setminus\mathcal{H})/\Gamma})^*$ due to Looijenga.  Roughly speaking, it involves blowing up the B-B boundary where $\mathcal{H}/\Gamma$ meets it, and then blowing down $\mathcal{H}/\Gamma$.

\subsubsection{The general case}

When $D/\Gamma$ is not algebraic, the obvious question is ``what about the image of the period map?''  While Griffiths and Sommese were able to settle the cases of $\M$ compact and $\phi(\M)$ smooth (respectively), the general problem remained open until 2018, when Bakker, Brunebarbe and Tsimerman \cite{BBT} used algebraization results in $o$-minimal geometry to show that $\phi(\M)$ is always quasi-projective.\footnote{The assumption in \cite{BBT} that $\Gamma$ is arithmetic is easily removed by appealing to Griffiths's Properness Theorem.}  They also show that the restriction of the Griffiths line bundle $L_{D/\Gamma}$ to $\phi(\M)$ is ample.  Ongoing work of Green, Griffiths, Laza and Robles aims to construct a projective compactification $\overline{\phi(\M)}^*$ to which $L_{D/\Gamma}$ has an ample extension, and which is also Hodge-theoretically modular in the sense of recording the $\mathrm{Gr}^W_{\bullet}(\text{LMHS})$ with respect to a good compactification $\overline{\M}\supset \M$.  So this would yield a generalization of the B-B compactification to the nonclassical case.

An earlier construction of Kato and Usui \cite{KU} extends toroidal compactifications to the general setting. 
The object $D_{\Sigma}/\Gamma$ which they associate to a $\Gamma$-compatible fan $\Sigma$ in $\fg_{\RR}$ (by ``gluing in'' the $\{B(\sigma)/\Gamma_{\sigma}\backslash\}_{\sigma\in \Sigma/\Gamma}$) is, outside the classical case, only a \emph{partial} compactification (as a ``log manifold'') with ``slits'' (imagine ``compactifying'' $\CC^2\setminus\{x=0\}$ by adding in just the point $\{(0,0)\}$).  However, if $\overline{\M}\supset \M$ is a good compactification, and every monodromy cone $\sigma_x$ ($x\in \overline{\M}\setminus \M$) of the period map $\phi$ is contained in a cone of $\Sigma$, then they nevertheless obtain a \emph{compactification} $\phi_{\Sigma}\colon \overline{\M}\to D_{\Sigma}/\Gamma$ of $\phi$ by its ``LMHS mod $\Gamma$'', with image a separated compact algebraic space.\footnote{The original hope in \cite{KU} was to have a \emph{complete} fan, such that $D_{\Sigma}/\Gamma$ would compactify images of \emph{every} period map $\phi$ into $D/\Gamma$; but by work of Usui's student Watanabe, these need not exist as soon as $\mathcal{W}$ has integral manifolds of dimension $>1$.  However, it still seems likely that, for each $\phi$, one can produce a suitable fan; and this is perfectly consistent with the classical case, since there (and only there) $\mathrm{Id}_{D/\Gamma}$ is a period map.}  So if one has a Torelli theorem for $\phi$, and can prove (say) that $\Gamma.\{\sigma_x\mid x\in \overline{\M}\setminus \M\}$ is a fan, then one obtains a Hodge-theoretically modular compactification of $\M$.  While we are only aware of \emph{generic} Torelli results in the nonclassical case so far (such as those of Donagi and Voisin for projective hypersurfaces), this appears to be a promising line of inquiry for future research.

\section{Compactifying the examples}
\label{sec:Geo+HT}

We are near the end of our tour.  But we did make a promise about new routes to old haunts, or something like that.
According to legend, a certain speakeasy in a certain city had two doors to two different streets: an entrance for the police, and an exit for the patrons (address `86').
Without `taking out' our readers, we shall finally try to describe how the different compactifications and methods come together for our examples from \S\ref{sec:ExampleModuliTheory}. 

\subsection{Beyond Picard curves}
Deligne and Mostow showed in 1986, with additional contributions by Doran in 2004, that there are isomorphisms between certain GIT completions of $\M_{0,n}$ and Baily-Borel compactifications of $(n-3)$-dimensional ball quotients.  From their theorems one obtains a collection of weights $\mathbf{w} \in \mathbb{Q}^n_{>0}$ and integers $m\in\mathbb{N}$ such that 
$(\mathbb{B}_{n-3}/\Gamma_{\mathbf{w}})^{*} \cong  (\mathbb{P}^1)^n/ \! \! /_{\mathbf{w}}\SL_2\times S_m$ 
(where $S_m$ is the $m^{\text{th}}$ symmetric group).  For a list of these cases see Tables  2 and 3 in \cite{GKS}.  The isomorphisms compactify period maps $\M_{0,n}/S_m\hookrightarrow \mathbb{B}_{n-3}/\Gamma_{\mathbf{w}}$ associated to cyclic covers of $\mathbb{P}^1$ branched in a manner dictated by the configuration of weighted points.  In each case, the cyclic automorphism of the covering curve has an eigenspace in $H^1$ with Hodge numbers $(1,n-3)$; and this is what fails for $n>12$.

In a complementary direction, using the stable pairs described in  \S\ref{sec:geocompact}, we obtain the Hassett moduli space $\overline{\M}_{0, \mathbf{w + \epsilon}}$ of $n$-pointed rational curves with weights $\mathbf{w}+\epsilon$. This compactification of $\M_{0,n}$ is a smooth projective variety, and it admits a morphism  to 
$(\mathbb{P}^1)^n/ \! /_{\mathbf{w}}\SL_2(\mathbb{C})$. We also have a unique toroidal compactification,
$\overline{\mathbb{B}_{n-3}/\Gamma_{\mathbf{w}}}^{\tor}$,
discussed in \S\ref{sec:HodgeTheory}. Recent work of the authors with L. Schaffler \cite{GKS} found that there is an isomorphism between $\overline{\M}_{0, \mathbf{w + \epsilon}}/S_m$ and the toroidal compactification.  Our route leads us to the 
following commutative diagram for $n \leq 12$:
\begin{align*}
\xymatrix{
\overline{\M}_{0,n} \ar[r]
&
\overline{\M}_{0, \mathbf{w + \epsilon}}/S_m
\ar[d] \ar[r]^{\cong}
&
\overline{\BB_{n-3} / \Gamma_{\mathbf{w}} }^{\tor}
\ar[d] 
\\
&
(\mathbb{P}^1)^n/ \! \! /_{\mathbf{w}}\SL_2\times S_m
\ar[r]^{\cong}
&
(\BB_{n-3} / \Gamma_{\mathbf{w}})^{*}
}    
\end{align*}
\subsection{Cubic surfaces}
In  \S\ref{sec:cubicsurfaces}, 
we learned that there are injective period maps 
$\M_{\ord} \hookrightarrow \BB_4 / \Gamma(\sqrt{-3})$
and
$\M \hookrightarrow \BB_4 / \Gamma$. 
Allcock, Carlson and Toledo showed in 2000 that these maps extend to isomorphisms between the respective GIT  and Baily-Borel compactifications. 

Along the MMP route, Hacking, Keel, and Tevelev considered the moduli space parametrizing pairs 
$(X,D)$ where $X$ is a smooth \emph{non-Eckardt} cubic surface and $D$
is the sum of its 27 lines $L_i$ (which then intersect transversely). 
The divisor is necessary because $K_X+D$ is an ample divisor but 
$K_X$ is not.  The moduli space of such pairs is open in $\M_{\ord}$, and we can compactify it by adding stable pairs.  Their work shows that the resulting compact moduli space  $\overline{\M}_{\ord}^{\ksba}$ is a smooth projective variety admitting a morphism to $\overline{\M}_{\ord}^{\GIT}$.

Another approach is to consider weighted pairs
$(X, \sum_{i=1}^{27} (\frac{1}{9} + \epsilon) L_i)$,
where the coefficient is the ``smallest" one for which
$K_X + \sum_{i=1}^{27} (\frac{1}{9} + \epsilon ) L_i$
is ample.
Let $\overline{\M}_{\tiny\text{ord}}^{\text{\tiny W}}$ be the KSBA moduli space parametrizing such weighted pairs and their stable degenerations.  The authors' work with L. Schaffler \cite{GKS} shows that its normalization 
$( \overline{\M}_{\tiny\text{ord}}^{\text{\tiny W}})^{\nu}$ is isomorphic to both the unique toroidal compactification of the ball quotient and a compactification constructed by Naruki in 1982 using Cayley's cross-ratios.

We arrive then at the following diagram.
$$
\xymatrix{
\overline{\M}_{\ord}^{\ksba}
\ar[r]
&
( \overline{\M}_{\ord}^{\mathrm{W}})^{\nu}
\ar[d] \ar[rr]^{\cong \hspace{1cm}}
&&
\overline{\BB_4 / 
\Gamma
\left( \sqrt{-3} \right)
}^{\tor}
\ar[d]
\\
&
\overline{\M}^{\GIT}_{\ord}  
\ar[rr]^{\cong \hspace{1cm}} 
\ar[d]^{/W(E_6)}
&&
(\BB_4 / \Gamma ( \sqrt{-3}))^{*}
\ar[d]^{/W(E_6)}
\\
\overline{\M}^{\KE}
\ar[r]^{\cong}
&
\overline{\M}^{\GIT}  \ar[rr]^{\cong  \hspace{1cm}}
&&
( \BB_4/ \Gamma )^{*} 
}
$$
\subsection{K3 surfaces}
We are entering a region where much less is known, so we will focus on the case of degree two.  For constructing a functorial compactification in terms of stable pairs $(X,D)$, we need a choice of an ample divisor $D$.  Which one shall we use? One option described by Laza is to use a divisor in the degree 2 linear system.  However, this increases the dimension of the moduli problem by two.

A natural choice which does not increase the dimension is available for polarized K3 surfaces:  the ramification divisor $R\subset X$ for the 2:1 map $X \to \PP^2$ branched over $C_6 \subset \mathbb{P}^2$. 
To describe the stable degenerations for the pairs $(X, \epsilon R)$ is the same as describing them for the pairs 
$\left( \mathbb{P}^2, (\frac{1}{2} + \epsilon) C_6 \right)$.
This last moduli space (of plane curves) was studied for arbitrary degree by Hacking in 2003 and partially described in the case of sextic plane curves.  In 2019, Alexeev, Engel, and Thompson completely described the moduli space $\overline{\M}^{\tiny{K3,\text{slc}}}_2$ of pairs $(X, \epsilon R)$. Furthermore, they showed that there exists a toroidal compactification (for the so-called Coxeter fan $\Sigma^{\text{cox}}$) mapping to the normalization  $(\overline{\M}^{\tiny\text{slc}}_2)^{\nu}$. 
Along with the work of Looijenga and Shah, which extends the period map from   $\mathrm{Bl}_{\omega}( \overline{\M}^{\tiny{K3,\GIT}}_2 )$ to the Baily-Borel compactification, this results in the diagram:

\begin{align*}
\xymatrix{
( \overline{\M}^{K3, \text{slc}}_2 )^{\nu} \ar [rd] & (\overline{D/ \Gamma_2} )_{\Sigma^{\text{cox}}} \ar [l] \ar  [d]
\\
\Bl_{\omega}( \overline{\M}^{K3,\GIT}_2 ) \ar [r] \ar [d] \ar [rd] & (D/\Gamma_2)^*
\\
\overline{\M}^{K3,\GIT}_2 \ar [r]^{\cong} & ((D\setminus \mathcal{H}_{\omega})/\Gamma_2)^*
}
\end{align*}

\subsection{Cubic 4-folds}
The GIT compactification  $\overline{\M}^{\mathrm{GIT}}$ was described by Laza in 2007. In particular, he showed that a cubic fourfold $Y$ with isolated singularities is GIT-stable iff $Y$ has at worst ADE singularities. Thus, there is an open $\M^{\mathrm{ADE}} \subset \overline{\M}^{\mathrm{GIT}}$ parametrizing them. The  period map described in  \S\ref{sec:Cubic4folds} extends (by Griffiths properness \cite{CMSP}) to a regular morphism $\M^{\mathrm{ADE}} \to D / \Gamma$ whose image is the complement of an  arrangement known as $\mathcal{H}_{\infty}$. Subsequent work of Laza shows that there is an isomorphism from the GIT quotient to the Looijenga compactification associated to $\mathcal{H}_{ \infty }$. 

In another direction, a compactification can also be obtained via $K$-stability because smooth cubic fourfolds are Fano varieties admitting a KE metric.  
Liu showed in 2020 that the moduli space $\overline{\M}^{\mathrm{KE}}$ is isomorphic to the GIT quotient.  Thus, we arrive via three different roads to the same compactification:
\begin{align*}
\overline{\M}^{\mathrm{KE}}
 \cong
\overline{\M}^{\mathrm{GIT}}
\cong
({(D\setminus\mathcal{H}_{\infty} )/\Gamma})^*.
\end{align*}

\subsection{Abelian varieties} 
\label{sec:abelian}
We will not do justice to our last stop. Yet, we include it because $\A_g$ admits several well-studied toroidal  compactifications. They depend on an admissible fan covering of the rational closure of the space of positive definite symmetric $g \times g$-matrices. There are three classical fans: the second Voronoi, the perfect cone, and the central cone decomposition. They lead to three compactifications of $\A_g$ which coincide for $g \leq 3$ but are different in general.  Among those, the second Voronoi compactification has a modular interpretation due to Alexeev:  it is the normalization of the main
irreducible component of the moduli space of stable semi-abelic pairs (see also subsequent work by Olsson in the context of log-geometry). 

For which choices of a fan does the Torelli map $\M_g \to \A_g$ extend to a regular map $\overline{\M}_g \to (\overline{\A}_g)_{\Sigma}$? For the 2nd Voronoi fan, a positive answer was given by Mumford and Namikawa in 1976. 
In 2011, Alexeev and Brunyate give a positive answer for the perfect cone, and a negative one for the central cone when $g \geq 9$. The reader wishing to learn more about these developments, and moduli of abelian varieties more generally, may consult \cite{abelian} and the references therein.

\end{document}